\documentclass[11pt]{article}
\usepackage{latexsym,amssymb}
\usepackage{amsmath}
\usepackage{slashbox}

\usepackage{multirow,eepic}

\makeatletter
\def\senbun#1(#2)#3({\@senbun(#2)(}
\def\@senbun(#1,#2)(#3,#4){%
   \@tempdima#1\p@ \advance\@tempdima#3\p@
   \divide\@tempdima\tw@
   \@tempdimb#2\p@ \advance\@tempdimb#4\p@
   \divide\@tempdimb\tw@
   \edef\@senbun@temp{\noexpand\qbezier(#1,#2)%
      (\strip@pt\@tempdima,\strip@pt\@tempdimb)(#3,#4)}%
   \@senbun@temp}
\makeatother

\newtheorem{Theorem}{\bf Theorem}[section]
\newtheorem{Lemma}{\bf Lemma}[section]
\newtheorem{Proposition}{\bf Proposition}[section]
\newtheorem{Corollary}{\bf Corollary}[section]
\newtheorem{Remark}{\bf Remark}[section]
\newtheorem{Example}{\bf Example}[section]
\newtheorem{Definition}{\bf Definition}[section]

\newenvironment{theorem}{\begin{Theorem}$\!\!\!$}{\end{Theorem}}
\newenvironment{lemma}{\begin{Lemma}$\!\!\!$}{\end{Lemma}}
\newenvironment{proposition}{\begin{Proposition}$\!\!\!$}{\end{Proposition}}
\newenvironment{corollary}{\begin{Corollary}$\!\!\!$}{\end{Corollary}}

\numberwithin{equation}{section}

\begin{document}

\title{Sharp decay estimates in Lorentz spaces \\for nonnegative Schr\"odinger heat semigroups}
\author{
\\
Norisuke Ioku\\
Graduate School of Science and Engineering, Ehime University\\
Matsuyama, Ehime 790-8577, Japan,\vspace{5pt}\\
\quad\\
Kazuhiro Ishige
\\
Mathematical Institute, Tohoku University\\
Aoba, Sendai 980-8578, Japan,\\
\quad\\
and\\
\quad\\
Eiji Yanagida
\\
Department of Mathematics, Tokyo Institute Technology\\
Meguro-ku, Tokyo 152-8551, Japan}
\date{}
\maketitle
\newpage
\section{Introduction} 
Let $H:=-\Delta+V$ be a Schr\"odinger operator on $L^2({\bf R}^N)$, 
where $N\ge 2$ and $V\in L^r_{{\rm loc}}({\bf R}^N)$ with $r>N/2$. 
Assume that the operator $H$ is nonnegative, that is, 
\begin{equation}
\label{eq:1.1}
\int_{{\bf R}^N}\left\{|\nabla\phi|^2+V\phi^2\right\}dx\ge 0
\quad\mbox{for all}\quad\phi\in C_0^\infty({\bf R}^N).
\end{equation}
The operator $H$ is said to be subcritical 
if, for any $W\in C_0^\infty({\bf R}^N)$, the operator $H_\epsilon:=-\Delta+V-\epsilon W$ is nonnegative on $L^2({\bf R}^N)$ 
for any sufficiently small $\epsilon>0$. 
This is equivalent to that 
there exists, for any $y\in{\bf R}^N$, 
a positive solution $G(x,y)$ of 
$$
(-\Delta+V(x))G(x,y)=\delta(x-y)\quad\mbox{in}\quad{\bf R}^N,
$$
where $\delta$ is the Dirac delta function. 
If $H$ is not subcritical, 
then the operator $H$ is said to be critical. 
Nonnegative Schr\"odinger operators have been studied by many mathematicians 
since the pioneering work due to Simon~\cite{S} 
(see also \cite{DS}--\cite{IK05}, \cite{MS}-\cite{P3}, \cite{S2}--\cite{Zao2}, and references therein),
and the large time behavior of $L^q$ norms of the Schr\"odinger heat semigroup $e^{-tH}$ depends on 
whether $H$ is subcritical or not and on the behavior of positive harmonic functions for $H$ at the space infinity. 

In this paper we focus on a nonnegative Schr\"odinger operator $H:=-\Delta+V$ with a radially symmetric potential 
$V=V(|x|)$ behaving like 
$$
V(r)=\omega r^{-2}(1+o(1))\quad\mbox{as}\quad r\to\infty,
$$
where 
$$
\omega>-\omega_*\qquad\mbox{and}\qquad \omega_*:=(N-2)^2/4,
$$ 
and study the Schr\"odinger heat semigroup $e^{-tH}$ in the Lorentz spaces.
More precisely, we assume the following: 
$$
(V)\qquad\qquad
\left\{
\begin{array}{ll}
({\rm i}) & \mbox{$V=V(r)\in C^1([0,\infty))$};\vspace{5pt}
\qquad\qquad\qquad\qquad\\
({\rm ii}) & \mbox{there exist constants $\omega>-\omega_*$ and $a>0$ such that}\vspace{7pt}\\
 & \qquad\qquad\qquad\quad
 V(r)=\omega r^{-2}+O(r^{-2-a})\quad\mbox{as $r\to\infty$;}\qquad\qquad\qquad\vspace{7pt}\\
({\rm iii}) & \displaystyle{\sup_{r>1}\,r^3|V'(r)|<\infty},
\end{array}
\right.
$$
and make the complete table of the decay rates of  
$$
\|e^{-tH}\|_{(L^{p,\sigma}\to L^{q,\theta})}
:=\sup\left\{\,\frac{\|e^{-tH}\phi\|_{L^{q,\theta}({\bf R}^N)}}{\|\phi\|_{L^{p,\sigma}({\bf R}^N)}}\,:\,
\phi\in L^{p,\sigma}({\bf R}^N)\setminus\{0\}\right\}
$$
as $t\to\infty$. 
Here $\|e^{-tH}\|_{(L^{p,\sigma}\to L^{q,\theta})}$ is 
the operator norm of $e^{-tH}$ 
from the Lorentz space $L^{p,\sigma}({\bf R}^N)$ to $L^{q,\theta}({\bf R}^N)$,
where 
$$
(p,q,\sigma,\theta)\in\Lambda:=
\left\{1\le p\le q\le\infty,\,\sigma,\,\theta\in[1,\infty]:
\begin{array}{lll}
\sigma=1 & \mbox{if}\quad p=1,\ \sigma=\infty & \mbox{if}\quad p=\infty\vspace{3pt}\\
\theta=1 & \mbox{if}\quad q=1,\ \theta=\infty & \mbox{if}\quad q=\infty\vspace{3pt}\\
\sigma\le\theta & \mbox{if}\quad p=q &
\end{array}
\right\}.
$$
This paper is an improvement and a generalization of our previous paper~\cite{IIY}, 
where the decay rates of the operator norms of 
$e^{-tH}$ in the $L^p$ spaces were discussed. 
\vspace{5pt}

We introduce some notation. 
For any sets $\Xi$ and $\Sigma$, let $f=f(\xi,\sigma)$ and $h=h(\xi,\sigma)$ 
be maps from $\Xi\times\Sigma$ to $(0,\infty)$. 
Then we say 
$$
f(\xi,\sigma)\preceq h(\xi,\sigma)
\quad\mbox{for all}\quad\xi\in\Xi
$$
if, for any $\sigma\in\Sigma$, 
there exists a positive constant $C$ such that $f(\xi,\sigma)\le Ch(\xi,\sigma)$ 
for all $\xi\in\Xi$. 
In addition, we say 
$$
f(\xi,\sigma)\asymp h(\xi,\sigma)
\quad\mbox{for all}\quad\xi\in\Xi
$$ 
if $f(\xi,\sigma)\preceq h(\xi,\sigma)$ and $f(\xi,\sigma)\succeq h(\xi,\sigma)$ 
for all $\xi\in\Xi$. 
Let $B(x,r):=\{y\in{\bf R}^N:\, |y-x|<r\}$  
for $x\in{\bf R}^N$ and $r>0$. 
For any measurable set $E\subset{\bf R}^N$, we denote by $\chi_E$ the characteristic function of $E$. 
\vspace{5pt}

Let $H:=-\Delta+V$ be a nonnegative operator on $L^2({\bf R}^N)$ and assume condition~$(V)$. 
Then there exists a positive radially symmetric harmonic function $U=U(|x|)$ for the operator $H$, that is, 
$$
U>0,\qquad -\Delta U+V(|x|)U=0\quad\mbox{in}\quad{\bf R}^N,
$$
and it satisfies 
\begin{equation}
\label{eq:1.2}
\lim_{r\to\infty}r^A U(r)=1,
\end{equation}
where 
\begin{equation}
\label{eq:1.3}
A:=
\left\{
\begin{array}{ll}
\displaystyle{\frac{N-2-\sqrt{(N-2)^2+4\omega}}{2}} & \quad\mbox{if $H$ is subcritical,}\vspace{8pt}\\
\displaystyle{\frac{N-2+\sqrt{(N-2)^2+4\omega}}{2}} & \quad\mbox{if $H$ is critical.}
\end{array}
\right.
\end{equation}
(See \cite[Theorem~5.7]{Murata03}.) 
Here we remark the following:  
\begin{itemize}
  \item \quad $\omega A\le 0$ and $A<N/2$ if $H$ is subcritical; 
  \item \quad$A>0$ if $H$ is critical; 
  \item \quad $U\not\in L^2({\bf R}^N)$ is equivalent to $A\le N/2$. 
\end{itemize}
For any $1\le p\le q\le\infty$, let 
$\|e^{-tH}\|_{(L^p\to L^q)}$ be the operator norm of the Schr\"odinger heat semigroup $e^{-tH}$ from $L^p({\bf R}^N)$ to $L^q({\bf R}^N)$, 
that is, 
$$
\|e^{-tH}\|_{(L^p\to L^q)}
:=\sup\left\{\,\frac{\|e^{-tH}\phi\|_{L^q({\bf R}^N)}}{\|\phi\|_{L^p({\bf R}^N)}}\,:\,
\phi\in L^p({\bf R}^N)\setminus\{0\}\right\}. 
$$
It follows form the nonnegativity of the operator $H$ that 
\begin{equation}
\label{eq:1.4}
\|e^{-tH}\|_{(L^2\to L^2)}\le 1,\qquad t>0. 
\end{equation}

Generally, 
the decay of the operator norms $\|e^{-tH}\|_{(L^p\to L^q)}$ as $t\to\infty$ depends 
on the behavior of the positive harmonic functions at the space infinity, 
and it has been studied in several papers 
(see e.g. \cite{DS}--\cite{IK05} and \cite{S}). 
Among others,  
the authors of~\cite{IIY} 
studied the decay rates of $\|e^{-tH}\|_{(L^p\to L^q)}$ as $t\to\infty$ 
under the assumption:  
$$
(\tilde{V})\qquad\qquad
\left\{
\begin{array}{ll}
({\rm i}) & \mbox{$V=V(r)\in C^1([0,\infty))$ and $V'\ge 0$, $V\le 0$ in $[0,\infty)$};\vspace{5pt}
\qquad\qquad\qquad\quad\\
({\rm ii}) & \mbox{there exist constants $\omega\in(-\omega_*,0]$ and $a>0$ such that}\vspace{7pt}\\
 & \qquad\qquad\qquad\quad
 V(r)=\omega r^{-2}+O(r^{-2-a})\quad\mbox{as $r\to\infty$;}\qquad\qquad\vspace{7pt}\\
({\rm iii}) & \displaystyle{\sup_{r>1}\,r^3|V'(r)|<\infty}.
\end{array}
\right.
$$
In particular, they gave the sharp decay rates of $\|e^{-tH}\|_{(L^p\to L^q)}$ as $t\to\infty$ 
for all $1\le p\le q\le\infty$ in the case $A<N/2$, and 
proved the following. 
\vspace{3pt}
\newline
(a)  Let $A=0$. Then 
$$
\|e^{-tH}\|_{(L^p\to L^q)}\asymp t^{-\frac{N}{2}(\frac{1}{p}-\frac{1}{q})}\asymp\|e^{t\Delta}\|_{(L^p\to L^q)}
$$
for all $t\ge 2$.
\vspace{3pt}
\newline
(b)  
Let $0<A<N/2$ and set 
$$
\alpha:=\frac{N}{N-A}\qquad\mbox{and}\qquad \beta:=\frac{N}{A}. 
$$
Then 
$$
\|e^{-tH}\|_{(L^p\to L^q)}\asymp\,\eta_{p,q}(t)
$$
for all $t\ge 2$, where $\eta_{p,q}(t)$ is defined by the following. 
\vspace{10pt}
\newline
\begingroup
\renewcommand{\arraystretch}{1.6}
\begin{table}[h]
\begin{tabular}{|c||c|c|c|c|c|}
\hline 
 $L^p\to L^q$& $1\le p<\alpha$ & $p=\alpha$ & $\alpha<p<\beta$ & $p=\beta$ & $\beta<p\le\infty$ 
 \\
\hline\hline $1\le q<\alpha$ & $t^{-\frac{N}{2q'}+\frac{A}{2}}$ 
& \senbun(87,15)(0,-6) 
& \senbun(85,15)(0,-6) 
& \senbun(50,15)(0,-6) 
& \senbun(64,15)(0,-6) 
\\
\hline $q=\alpha$ & $t^{-\frac{N}{2q'}+\frac{A}{2}}$ & $(\log t)^{\frac{A}{N}}$ 
& \senbun(85,15)(0,-6) 
& \senbun(50,15)(0,-6) 
& \senbun(64,15)(0,-6) 
\\
\hline $\alpha<q<\beta$ & $t^{-\frac{N}{2q'}+\frac{A}{2}}$ & $t^{-\frac{N}{2q'}+\frac{A}{2}}(\log t)^{\frac{A}{N}}$ & $t^{-\frac{N}{2}(\frac{1}{p}-\frac{1}{q})}$ 
& \senbun(50,15)(0,-6) 
& \senbun(64,15)(0,-6) 
\\
\hline $q=\beta$ & $t^{-\frac{N}{2}+A}(\log t)^{\frac{A}{N}}$ & $t^{-\frac{N}{2}+A}(\log t)^{\frac{2A}{N}}$ & 
$t^{-\frac{N}{2p}+\frac{A}{2}}(\log t)^{\frac{A}{N}}$ & $(\log t)^{\frac{A}{N}}$ 
& \senbun(64,15)(0,-6) 
\\
\hline $\beta<q\le\infty$ & $t^{-\frac{N}{2}+A}$ & $t^{-\frac{N}{2}+A}(\log t)^{\frac{A}{N}}$ 
& $t^{-\frac{N}{2p}+\frac{A}{2}}$ & $t^{-\frac{N}{2p}+\frac{A}{2}}$ & $t^{-\frac{N}{2p}+\frac{A}{2}}$ 
\\
\hline 
\end{tabular}
\caption{Decay rates of $\|e^{-tH}\|_{(L^p\to L^q)}$}
\end{table}
\endgroup
\vspace{7pt}

\noindent
Here, for any $q \in[1,\infty]$, 
we denote by $q'$ the H\"older conjugate number of $q$, that is, 
$$
q':=\frac{q}{q-1}\quad\mbox{if}\quad q\in(1,\infty),
\qquad 
q':=1\quad\mbox{if}\quad q=\infty,
\qquad 
q':=\infty\quad\mbox{if}\quad q=1.
$$ 
Furthermore, we remark that $1<\alpha<2<\beta$ and $\beta=\alpha'$.  
\vspace{5pt}

In this paper 
we eliminate the restriction of the sign of $V$ and $V'$ from condition~$(\tilde{V})$, 
and give the sharp decay rates of the operator norms of $e^{-tH}$ in the Lorentz spaces, 
which are more general function spaces than the $L^p$ spaces. 
In particular, we prove that, 
for the case where $A>0$ and 
$p=\alpha$ or $q=\beta$,
the decay rates of $\|e^{-tH}\|_{(L^{p,\sigma}\to L^{q,\theta})}$ 
depend on the second exponents $\sigma$ and $\theta$ of the Lorentz spaces $L^{p,\sigma}({\bf R}^N)$ and $L^{q,\theta}({\bf R}^N)$. 
As far as we know, there are no results pointing out 
the importance of the second exponents of the Lorentz spaces  
in the study of the behavior of the Schr\"odinger heat semigroups. 
\vspace{3pt}

Now we are ready to state the main result of this paper. 
We remark that $L^{p,p}({\bf R}^N)=L^p({\bf R}^N)$ for $1\le p\le \infty$ (see \eqref{eq:2.2}). 
\begin{theorem}
\label{Theorem:1.1} 
Let $N\ge 2$ and $H:=-\Delta+V$ be a nonnegative Schr\"odinger operator on $L^2({\bf R}^N)$. 
Assume condition~$(V)$ and $A<N/2$. 
Then, for any $(p,q,\sigma,\theta)\in\Lambda$, 
the following holds. 
\begin{itemize}
  \item[{\rm(I)}] 
	Let $A\le 0$. Then 
	$$
	\|e^{-tH}\|_{(L^{p,\sigma}\to L^{q,\theta})}\asymp t^{-\frac{N}{2}(\frac{1}{p}-\frac{1}{q})}
	$$
	for all $t\ge 2$. 
  \item[{\rm(II)}]  
  	Let $A>0$. 
	\begin{itemize}
		\item[{\rm (i)}] If $1\le p<\alpha$, then 
		$$
		\|e^{-tH}\|_{(L^{p,\sigma}\to L^{q,\theta})}\asymp
		\left\{
		\begin{array}{ll}
		t^{-\frac{N}{2}(1-\frac{1}{q})+\frac{A}{2}} & \mbox{if}\quad p\le q<\beta,\\
		t^{-\frac{N}{2}+A}(\log t)^{\frac{1}{\theta}} & \mbox{if}\quad q=\beta,\\
		t^{-\frac{N}{2}+A} & \mbox{if}\quad \beta<q\le\infty,
		\end{array}
		\right.
		$$
		for all $t\ge 2$.\vspace{3pt}
		\item[{\rm (ii)}] If $p=\alpha$, then
		$$
		\|e^{-tH}\|_{(L^{p,\sigma}\to L^{q,\theta})}\asymp
		\left\{
		\begin{array}{ll}
		t^{-\frac{N}{2}(1-\frac{1}{q})+\frac{A}{2}}(\log t)^{\frac{1}{\sigma'}} & \mbox{if}\quad \alpha \le q<\beta,\vspace{3pt}\\
		t^{-\frac{N}{2}+A}(\log t)^{\frac{1}{\theta}+\frac{1}{\sigma'}} & \mbox{if}\quad q=\beta,\vspace{3pt}\\
		t^{-\frac{N}{2}+A}(\log t)^{\frac{1}{\sigma'}} & \mbox{if}\quad \beta<q\le\infty,
		\end{array}
		\right.
		$$
		for all $t\ge 2$.\vspace{3pt}
		\item[{\rm (iii)}] If $\alpha<p<\beta$, then 
		$$
		\|e^{-tH}\|_{(L^{p,\sigma}\to L^{q,\theta})}\asymp
		\left\{
		\begin{array}{ll}
		t^{-\frac{N}{2}(\frac{1}{p}-\frac{1}{q})} & \mbox{if}\quad p\le q<\beta,\\
		t^{-\frac{N}{2p}+\frac{A}{2}}(\log t)^{\frac{1}{\theta}} & \mbox{if}\quad q=\beta,\\
		t^{-\frac{N}{2p}+\frac{A}{2}} & \mbox{if}\quad \beta<q\le\infty,
		\end{array}
		\right.
		$$
		for all $t\ge 2$.\vspace{3pt}
		\item[{\rm (iv)}] If $p=\beta$, then 
		$$
		\|e^{-tH}\|_{(L^{p,\sigma}\to L^{q,\theta})}\asymp
		\left\{
		\begin{array}{ll}
		(\log t)^{\frac{1}{\theta}} & \mbox{if}\quad q=\beta,\vspace{3pt}\\
		t^{-\frac{N}{2p}+\frac{A}{2}}	 & \mbox{if}\quad \beta<q\le \infty,	
		\end{array}
		\right.
		$$
		for all $t\ge 2$.\vspace{3pt}
		\item[{\rm (v)}] If $\beta<p\le \infty$, then 
		$$
		\|e^{-tH}\|_{(L^{p,\sigma}\to L^{q,\theta})}\asymp t^{-\frac{N}{2p}+\frac{A}{2}}
		$$
		for all $t\ge 2$. 
		\end{itemize}
\end{itemize}
\end{theorem}
By Theorem~\ref{Theorem:1.1} 
we have the following table on the decay rates of $\|e^{-tH}\|_{(L^{p,\sigma}\to L^{q,\theta})}$ in the case $0<A<N/2$. 
\vspace{8pt}

\noindent
\begingroup
\renewcommand{\arraystretch}{1.6}
\begin{table}[!h]
\begin{tabular}{|c||c|c|c|c|c|}
\hline 
$L^{p,\sigma}\to L^{q,\theta}$ & $1\le p<\alpha$ & $p=\alpha$ & $\alpha<p<\beta$ & $p=\beta$ & $\beta<p\le \infty$ 
 \\
\hline\hline $1\le q<\alpha$ & $t^{-\frac{N}{2q'}+\frac{A}{2}}$ 
& \senbun(96,15)(0,-6) 
& \senbun(82,15)(0,-6) 
& \senbun(47,15)(0,-6) 
& \senbun(64,15)(0,-6) 
\\
\hline $q=\alpha$ & $t^{-\frac{N}{2q'}+\frac{A}{2}}$ & $(\log t)^{\frac{1}{\sigma'}}$ 
& \senbun(82,15)(0,-6) 
& \senbun(47,15)(0,-6) 
& \senbun(64,15)(0,-6) 
\\
\hline $\alpha<q<\beta$ & $t^{-\frac{N}{2q'}+\frac{A}{2}}$ & $t^{-\frac{N}{2q'}+\frac{A}{2}}(\log t)^{\frac{1}{\sigma'}}$ & $t^{-\frac{N}{2}(\frac{1}{p}-\frac{1}{q})}$ 
& \senbun(47,15)(0,-6) 
& \senbun(64,15)(0,-6) 
\\
\hline $q=\beta$ & $t^{-\frac{N}{2}+A}(\log t)^{\frac{1}{\theta}}$ & $t^{-\frac{N}{2}+A}(\log t)^{\frac{1}{\theta}+\frac{1}{\sigma'}}$ & 
$t^{-\frac{N}{2p}+\frac{A}{2}}(\log t)^{\frac{1}{\theta}}$ & $(\log t)^{\frac{1}{\theta}}$ 
& \senbun(64,15)(0,-6) 
\\
\hline $\beta<q\le\infty$ & $t^{-\frac{N}{2}+A}$ & $t^{-\frac{N}{2}+A}(\log t)^{\frac{1}{\sigma'}}$ 
& $t^{-\frac{N}{2p}+\frac{A}{2}}$ & $t^{-\frac{N}{2p}+\frac{A}{2}}$ & $t^{-\frac{N}{2p}+\frac{A}{2}}$ 
\\
\hline 
\end{tabular}
\caption{Decay rates of $\|e^{-tH}\|_{(L^{p,\sigma}\to L^{q,\theta})}$}
\end{table}
\endgroup
\vspace{5pt}

\noindent
Furthermore, as a corollary of Theorem~\ref{Theorem:1.1}, we have: 
\begin{corollary}
\label{Corollary:1.1}
Let $N\ge 2$ 
and $H:=-\Delta+V$ be a nonnegative Schr\"odinger operator on $L^2({\bf R}^N)$.  
Assume condition~$(V)$ and $A<N/2$. 
Then 
$$
\|e^{-tH}\|_{(L^p\to L^q)}\asymp t^{-\frac{N}{2}(\frac{1}{p}-\frac{1}{q})},
\qquad t\ge 2,
$$
for all $1\le p\le q\le\infty$ 
if and only if $H$ is subcritical and $\omega\ge 0$. 
\end{corollary}
Corollary~\ref{Corollary:1.1} immediately follows from Theorem~\ref{Theorem:1.1}. 
\vspace{5pt}

We explain the idea of the proof of Theorem~\ref{Theorem:1.1}. 
Let $\delta>0$ and define 
\begin{equation}
\label{eq:1.5}
\chi_\delta(x,t):=0\quad\mbox{if}\quad |x|\le\delta(1+t)^{1/2},
\qquad
\chi_\delta(x,t):=1\quad\mbox{if}\quad |x|>\delta(1+t)^{1/2}. 
\end{equation}
We construct a supersolution of 
\begin{equation}
\label{eq:1.6}
\partial_t u=\Delta u-V(|x|)u\quad\mbox{in}\quad{\bf R}^N\times(0,\infty)
\end{equation}
with $u(x,0)=1$ in ${\bf R}^N$, 
and prove that  
$$
\|\chi_\delta(t)e^{-tH}\|_{(L^\infty\to L^\infty)}\le C,\qquad t\ge 2,
$$
for some constant $C$. 
This together with \eqref{eq:1.4} and a Marcinkiewicz type interpolation theorem in the Lorentz spaces 
implies that 
\begin{equation}
\label{eq:1.7}
\|\chi_\delta(t)e^{-tH}\|_{(L^{p,\sigma}\to L^{p,\sigma})}\le C',\qquad t\ge 2,
\end{equation}
for some constant $C'$. 
Furthermore, applying the $L^\infty_{loc}$ estimates for 
parabolic equations and using another supersolution of \eqref{eq:1.6}, 
we obtain the upper decay estimates of $\|e^{-tH}\|_{(L^{p,\sigma}\to L^{q,\theta})}$. 
On the other hand, 
the lower decay estimates of $\|e^{-tH}\|_{(L^{p,\sigma}\to L^{q,\theta})}$ 
are obtained by modification of the arguments in \cite{IIY} and \cite{IK05}. 

In our previous paper~\cite{IIY},  
we assumed that $V=V(r)$ is nonpositive and monotone increasing in $[0,\infty)$ (see condition~$(\tilde{V})$), 
and proved the inequality 
\begin{equation}
\label{eq:1.8}
\|e^{-tH}\phi\|_{L^{p,\infty}({\bf R}^N)}\le\|e^{-tH}\phi^\sharp\|_{L^{p,\infty}({\bf R}^N)},\qquad t>0,
\end{equation}
with the aid of \cite{ALT}. 
Here $\phi^\sharp$ is the spherical rearrangement of $\phi$ (see Section~2). 
The inequality~\eqref{eq:1.8} is a crucial ingredient in \cite{IIY} and 
its proof in \cite{IIY} requires the restriction of the sign of $V$ and $V'$. 
In this paper, without the use of the inequality~\eqref{eq:1.8}, 
we study the decay rates of  $\|e^{-tH}\|_{(L^{p,\sigma}\to L^{p,\sigma})}$.  
This enables us to obtain 
the sharp decay rates of $\|e^{-tH}\|_{(L^{p,\sigma}\to L^{q,\theta})}$ 
under condition~$(V)$, which is weaker than condition~$(\tilde{V})$. 
\vspace{3pt}

The rest of this paper is organized as follows. 
In Section~2 we recall some properties of the Lorentz spaces and 
some preliminary results on the Schr\"odinger operator $H$. 
In Sections~3 and 4 we give decay estimates of $\|e^{-tH}\|_{(L^{p,\sigma}\to L^{q,\theta})}$ 
by using supersolutions of \eqref{eq:1.6} with the aid of $L^\infty_{loc}$ estimates for 
parabolic equations 
and a Marcinkiewicz type interpolation theorem in the Lorentz spaces. 
In Section~5 we give lower estimates of $\|e^{-tH}\|_{(L^{p,\sigma}\to L^{q,\theta})}$, 
and complete the proof of Theorem~\ref{Theorem:1.1}. 
\section{Preliminaries} 
In this section we recall some properties of  
the Lorentz spaces and nonnegative Schr\"odinger operators. 
\vspace{3pt}

For any measurable function $\phi$ in ${\bf R}^N$, 
we denote by $\mu=\mu(\lambda)$ the distribution function of $\phi$, that is, 
$$
\mu(\lambda):=\left|\{x\,:\,|\phi(x)|>\lambda\}\right| 
\qquad (\lambda > 0).
$$
We define 
the non-increasing rearrangement $\phi^*$ of $\phi$
and the spherical rearrangement $\phi^{\sharp}$ of $\phi$ by 
$$
\phi^{*}(s):=\inf\{\lambda>0\,:\,\mu(\lambda)\le s\},
\qquad
\phi^\sharp(x):=\phi^*(c_N|x|^N),
$$
for $s>0$ and $x\in {\bf R}^N$, respectively,
where $c_N$ is the volume of the unit ball in ${\bf R}^N$. 
Then, for any $1\le p\le\infty$ and $1\le \sigma\le \infty$, 
we define the Lorentz space $L^{p,\sigma}({\bf R}^N)$ by 
$$
L^{p,\sigma}({\bf R}^N):=\{\phi\,:\, \mbox{$\phi$ is measurable on ${\bf R}^N$},\,\,\, \|\phi\|_{L^{p,\sigma}}<\infty\},
$$
where 
\begin{equation}
\label{eq:2.1}
\|\phi\|_{L^{p,\sigma}}:=
 \left\{
\begin{array}{ll}
  \displaystyle{\biggr(
   \int_{{\bf R}^N}\left(|x|^{N/p}\phi^{\sharp}(x)\right)^{\sigma}\frac{dx}{|x|^N}
  \biggr)^{1/\sigma}}\quad & \mbox{if}\quad 1\le \sigma<\infty, \vspace{5pt}\\
 \displaystyle{\sup_{x\in {\bf R}^N}\,|x|^{N/p}\phi^{\sharp}(x)}\quad & \mbox{if}\quad\sigma=\infty.
\end{array}
 \right.
\end{equation}
The Lorentz spaces have the following properties:  
\begin{eqnarray}
\label{eq:2.2}
  & & L^{p,p}({\bf R}^N)=L^p({\bf R}^N)\mbox{ if $1\le p\le \infty$};\vspace{3pt}\\ 
\label{eq:2.3}
  & & L^{p,\sigma}({\bf R}^N)\subset L^{p,\rho}({\bf R}^N)\mbox{ if $1\le p<\infty$ and $1\le \sigma\le \rho \le \infty$};\vspace{3pt}\\ 
\label{eq:2.4}
  & & L^{p,\sigma}({\bf R}^N)'=L^{p',\sigma'}({\bf R}^N)\mbox{ if $(p,p,\sigma,\sigma)\in\Lambda$}.
\end{eqnarray}
Here $L^{p,\sigma}({\bf R}^N)'$ is the associate space of $L^{p,\sigma}({\bf R}^N)$. 
See e.g. \cite[Theorem~4.7, Chapter~4]{BS}.
\vspace{7pt}

We state a Marcinkiewicz type interpolation theorem in the Lorentz spaces. 
Proposition~\ref{Proposition:2.1} follows from \cite[Theorem~1.12, Chapter~5]{BS} and \cite[Theorem~5.3.1]{BL}. 
\begin{proposition}
\label{Proposition:2.1}
Let $(p_0,q_0,\sigma_0,\theta_0)\in \Lambda$ and $(p_1,q_1,\sigma_1,\theta_1)\in \Lambda$.
For $0<\eta<1$, set 
$$
\frac{1}{p}:=\frac{1-\eta}{p_0}+\frac{\eta}{p_1}, 
\quad
\frac{1}{q}:=\frac{1-\eta}{q_0}+\frac{\eta}{q_1}.
$$
Let $T$ be a bounded linear operator from 
$L^{p_0,\sigma_0}({\bf R}^N)$ to $L^{q_0,\theta_0}({\bf R}^N)$ and
from $L^{p_1,\sigma_1}({\bf R}^N)$ to $L^{q_1,\theta_1}({\bf R}^N)$, and define 
$$
M_0:=\sup_{f\in L^{p_0,\sigma_0}({\bf R}^N)\setminus\{0\}}
\frac{\|Tf\|_{L^{q_0,\theta_0}({\bf R}^N)}}{\|f\|_{L^{p_0,\sigma_0}({\bf R}^N)}},
\qquad
M_1:=\sup_{f\in L^{p_1,\sigma_1}({\bf R}^N)\setminus\{0\}}
\frac{\|Tf\|_{L^{q_1,\theta_1}({\bf R}^N)}}{\|f\|_{L^{p_1,\sigma_1}({\bf R}^N)}}.
$$
Then the following holds. 
\begin{itemize}
\item[\rm{(i)}] 
If $p_0\neq p_1$ and $q_0\neq q_1$, then 
$$ 
\|Tf\|_{L^{q,\sigma}}\le M_0^{1-\eta}M_1^{\eta} \|f\|_{L^{p,\sigma}}\quad 
\mbox{for every\quad  $1\le \sigma \le \infty$}.
$$
\item[\rm{(ii)}] 
If $p_0\neq p_1$, $q_0= q_1$, and $\theta_0=\theta_1=\theta$, then 
$$
\|Tf\|_{L^{q,\theta}}\le M_0^{1-\eta}M_1^{\eta} \|f\|_{L^{p,\sigma}}\quad 
\mbox{for every\quad  $1\le \sigma \le \infty$}. 
$$
\item[\rm{(iii)}] 
If $p_0= p_1$, $q_0= q_1$, and $\theta_0=\theta_1=\theta$, then 
$$
\|Tf\|_{L^{q,\theta}}\le M_0^{1-\eta}M_1^{\eta} \|f\|_{L^{p,\sigma}}\quad 
\mbox{for\quad $\dfrac{1}{\sigma}=\dfrac{1-\eta}{\sigma_0}+\dfrac{\eta}{\sigma_1}$}.
$$
\end{itemize}
\end{proposition}

We prove the following proposition on the Schr\"odinger heat semigroup $e^{-tH}$. 
\begin{proposition}
\label{Proposition:2.2}
Assume the same conditions as in Theorem~{\rm\ref{Theorem:1.1}}. 
Let $(p,q,\sigma,\theta)\in\Lambda$. 
Then, for any $T>0$, 
\begin{equation}
\label{eq:2.5}
\|e^{-tH}\|_{(L^{p,\sigma}\to L^{q,\theta})}\preceq t^{-\frac{N}{2}(\frac{1}{p}-\frac{1}{q})}
\end{equation}
for all $t\in(0,T)$. 
Furthermore, 
\begin{equation}
\label{eq:2.6}
\|e^{-tH}\|_{(L^{p,\sigma}\to L^{q,\theta})}
\asymp
\|e^{-tH}\|_{(L^{q',\theta'}\to L^{p', \sigma'})}
\end{equation}
for all $t>0$, where 
$p'$, $q'$, $\sigma'$, and $\theta'$ are the H\"older conjugate numbers of $p$, $q$, $\sigma$, and $\theta$, respectively. 
\end{proposition}
{\bf Proof.} 
Let $T>0$ and $\phi\in L^{p,\sigma}({\bf R}^N)$ with $\phi\not\equiv 0$ in ${\bf R}^N$. 
Due to  condition~$(V)$, we see that $V\in L^\infty({\bf R}^N)$, 
and we can define the function $v$ by 
$$
v(x,t):=e^{t\|V\|_{L^\infty({\bf R}^N)}}[e^{t\Delta}|\phi|](x). 
$$
Since 
$v$ satisfies 
$$
\partial_t v=\Delta v+\|V\|_{L^\infty({\bf R}^N)}v\quad\mbox{in}\quad{\bf R}^N\times(0,\infty),
\qquad
v(x,0)=|\phi(x)|\quad\mbox{in}\quad{\bf R}^N,
$$
by the comparison principle we have 
$$
\left|[e^{-tH}\phi](x)\right|\le v(x,t)=e^{t\|V\|_{L^\infty({\bf R}^N)}}[e^{t\Delta}|\phi|](x)
$$
for all $(x,t)\in{\bf R}^N\times(0,\infty)$. 
This implies
$$
\|e^{-tH}\|_{(L^{p,\sigma}\to L^{q,\theta})}
=\sup_{\phi\in L^{p,\sigma}({\bf R}^N)\setminus\{0\}}
\frac{\|e^{-tH}\phi\|_{L^{q,\theta}}}{\|\phi\|_{L^{p,\sigma}}}
\preceq \|e^{t\Delta}\|_{(L^{p,\sigma}\to L^{q,\theta})}
\asymp t^{-\frac{N}{2}(\frac{1}{p}-\frac{1}{q})}
$$
for all $t\in(0,T)$, and we have \eqref{eq:2.5}. 
On the other hand, similarly to \cite[Proposition~2.2]{IIY}, it follows that  
$$
\int_{{\bf R}^N}[e^{-tH}\phi](x)\psi(x)dx=\int_{{\bf R}^N}\phi[e^{-tH}\psi](x)dx,\qquad t>0,
$$
for all $\phi\in L^{p,\sigma}({\bf R}^N)$ and $\psi\in L^{q',\theta'}({\bf R}^N)$. 
This together with \eqref{eq:2.4} 
implies \eqref{eq:2.6}, and the proof
 is complete.
$\Box$
\vspace{5pt}

Proposition~\ref{Proposition:2.3} is concerned with the behavior of positive harmonic functions for the operator $H$.
\begin{proposition}
\label{Proposition:2.3}
Assume the same conditions as in Theorem~{\rm\ref{Theorem:1.1}}. 
Then there exists a radially symmetric positive function $U=U(|x|)$ in ${\bf R}^N$ such that 
\begin{eqnarray}
\label{eq:2.7}
 & & \Delta U-V(|x|)U=0\quad\mbox{in}\quad{\bf R}^N,\vspace{3pt}\\
 \label{eq:2.8}
 & & U(r)=r^{-A}(1+o(1))\quad\mbox{as}\quad r\to\infty,\vspace{3pt}\\
\label{eq:2.9}
 & & U'(r)=-Ar^{-A-1}(1+o(1))\quad\mbox{as}\quad r\to\infty,
\end{eqnarray}
where $A$ is the constant given in \eqref{eq:1.3}. 
In particular, 
\begin{equation}
\label{eq:2.10}
U(r)\asymp(1+r)^{-A},\qquad r\ge 0.
\end{equation}
\end{proposition}
{\bf Proof.} 
Due to condition~$(V)$, 
it follows from \cite[Theorem~5.7]{Murata03} that 
there exists a radially symmetric positive function $U$ satisfying \eqref{eq:2.7} and \eqref{eq:2.8}. 
Furthermore, by \eqref{eq:2.8} and the positivity of $U$ we have \eqref{eq:2.10}. 
Moreover, by a similar argument as in the proof of (1.15) in \cite[Theorem~1.1]{IK04} 
we obtain \eqref{eq:2.9}. Thus Proposition~\ref{Proposition:2.3} follows. 
$\Box$
\vspace{5pt}

At the end of this section,
we state a proposition on supersolutions of \eqref{eq:1.6}. 
\begin{proposition}
\label{Proposition:2.4}
Assume the same conditions as in Theorem~{\rm\ref{Theorem:1.1}}. 
For any $\gamma_1$ and $\gamma_2\in{\bf R}$, take a constant $c>1$ such that 
$$
\zeta(t):=
(1+t)^{\gamma_1+\frac{A}{2}}[\log(c+t)]^{\gamma_2}
$$
is monotone in $[0,\infty)$. 
Then, for any $T>0$ and any sufficiently small $\epsilon>0$,  
there exist a constant $C$ and  a function $w(x,t)$ such that 
\begin{eqnarray*}
 & & \partial_t w\ge\Delta w-V(|x|)w
\quad\mbox{in}\quad{\bf R}^N\times(0,\infty),\\
 & & 0<w(x,t)\le C\zeta(t)U(|x|)
\quad\mbox{in}\quad D_\epsilon(T),\\
 & & w(x,t)\ge (1+t)^{\gamma_1}[\log(2+t)]^{\gamma_2}
\quad\mbox{on}\quad \Gamma_\epsilon(T),
\end{eqnarray*}
where 
\begin{eqnarray*}
D_\epsilon(T) & := & 
\left\{(x,t)\in{\bf R}^N\times(T,\infty)\,:\,|x|<\epsilon(1+t)^{1/2}\right\},\\
\Gamma_\epsilon(T) & := & 
\left\{(x,t)\in{\bf R}^N\times(T,\infty)\,:\,|x|=\epsilon(1+t)^{1/2}\right\}\\
&\, &\hspace{25mm}\cup\left\{(x,T)\in{\bf R}^N\times\{T\}\,:\, |x|<\epsilon(1+T)^{1/2}\right\}. 
\end{eqnarray*}
\end{proposition}
Proposition~\ref{Proposition:2.4} is proved 
by the same argument as in the proof of \cite[Lemma~3.1]{IK04}.
\section{Decay estimates of $\|e^{-tH}\|_{(L^{p,\sigma}\to L^{p,\sigma})}$} 
%
This section is devoted to the proof of the following proposition, 
which gives the decay estimates of 
$\|e^{-tH}\|_{(L^{p,\sigma}\to L^{p,\sigma})}$ as $t\to\infty$. 
\begin{proposition}
\label{Proposition:3.1}
Assume the same conditions as in Theorem~{\rm\ref{Theorem:1.1}}. 
Let $(p,p,\sigma,\sigma)\in\Lambda$. 
\begin{itemize}
  \item[{\rm(I)}]
  Let $A\le 0$. Then 
 $$ 
 \|e^{-tH}\|_{(L^{p,\sigma}\to L^{p,\sigma})}\preceq 1
 $$
 for all $t\ge 2$. 
  \item[{\rm(II)}] 
  Let $A>0$. Then 
  $$
 \|e^{-tH}\|_{(L^{p,\sigma}\to L^{p,\sigma})}
 \preceq\left\{
 \begin{array}{ll}
 t^{-\frac{N}{2p'}+\frac{A}{2}} & \mbox{if}\quad 1\le p<\alpha,\vspace{3pt}\\
 (\log t)^{\frac{1}{\sigma'}} & \mbox{if}\quad p=\alpha,\vspace{3pt}\\
 1 & \mbox{if}\quad \alpha < p<\beta,\vspace{3pt}\\
 (\log t)^{\frac{1}{\sigma}} & \mbox{if}\quad p=\beta,\vspace{3pt}\\
 t^{-\frac{N}{2p}+\frac{A}{2}} & \mbox{if}\quad \beta<p\le\infty,
 \end{array}
 \right.
 $$
for all $t\ge 2$. 
\end{itemize}
\end{proposition}
In order to prove Proposition~\ref{Proposition:3.1}, 
we first prove the following lemma by using the comparison principle and 
the Marcinkiewicz interpolation theorem. 
\begin{lemma}
\label{Lemma:3.1}
Assume the same conditions as in Theorem~{\rm\ref{Theorem:1.1}}. 
Let $(p,p,\sigma,\sigma)\in\Lambda$ with $2<p\le \infty$. 
Then, for any $\delta>0$, 
\begin{equation}
\label{eq:3.1}
\|\chi_\delta(t)e^{-tH}\|_{(L^{p,\sigma}\to L^{p,\sigma})}\preceq 1
\end{equation}
for all $t>0$, where $\chi_\delta=\chi_\delta(x,t)$ is the function given in \eqref{eq:1.5}. 
\end{lemma}
{\bf Proof.} 
We consider the case $A>0$. 
By the same argument as in \cite[Lemma~3.2]{IIY} 
we can construct a supersolution $W_1=W_1(|x|,t)$ of \eqref{eq:1.6} satisfying 
$$
W_1(|x|,t)\asymp  
\left\{
\begin{array}{ll}
(T_1+t)^{\frac{A}{2}}U(|x|) & \mbox{if}\quad |x|\le R_1(T_1+t)^{1/2},\medskip\\
1 & \mbox{if}\quad |x|\ge R_1(T_1+t)^{1/2},
\end{array}
\right.\\
$$
for all $(x,t)\in{\bf R}^N\times(0,\infty)$, 
where $T_1$ and $R_1$ are some positive constants. 
Then, by \eqref{eq:2.10} we have 
$$
\begin{array}{ll}
\qquad\quad W_1(x,0)\succeq 1,\qquad & x\in{\bf R}^N,\vspace{5pt}\\
\chi_\delta(x,t)W_1(x,t)\preceq 1,\qquad & (x,t)\in{\bf R}^N\times(0,\infty),
\end{array}
$$
for any $\delta>0$. 
These yield 
\begin{eqnarray}
\|\chi_\delta(t)e^{-tH}\|_{(L^\infty\to L^\infty)}
 \!\!\!& \le &\!\!\! \|\chi_\delta(t)e^{-tH}1\|_{L^\infty({\bf R}^N)}
\preceq\|\chi_\delta(t)e^{-tH}W_1(0)\|_{L^\infty({\bf R}^N)}\nonumber\\
\label{eq:3.2}
 \!\!\!& \le &\!\!\! \|\chi_\delta(t)W_1(t)\|_{L^\infty({\bf R}^N)}
\preceq 1
\end{eqnarray}
for all $t>0$. 
On the other hand, since $H$ is nonnegative, we have 
\begin{equation}
\label{eq:3.3}
\|\chi_\delta(t)e^{-tH}\|_{(L^2\to L^2)}
\le\|e^{-tH}\|_{(L^2\to L^2)}\le 1
\end{equation}
for all $t>0$ (see also \eqref{eq:1.4}). 
Therefore, by \eqref{eq:3.2} and \eqref{eq:3.3} 
we apply Proposition~\ref{Proposition:2.1}~(i), 
and obtain \eqref{eq:3.1} in the case $A>0$. 
In the case $A=0$, since 
$U(|x|)\asymp 1$ in ${\bf R}^N$, 
taking $W_1(x,t)=U(|x|)$, 
we apply the same argument as in the case $A>0$ to obtain \eqref{eq:3.1}. 

It remains to prove \eqref{eq:3.1} in the case $A<0$.
Assume $A<0$. Let $T_2$ be a sufficiently large constant to be chosen later. 
Define 
\begin{eqnarray}
 & & v(|x|,t):=(T_2+t)^{\frac{A}{2}}\left[U(|x|)+\frac{A}{2}(T_2+t)^{-1}F(|x|)\right],\label{eq:3.4}\\
 & & F(|x|):=U(|x|)\int_0^{|x|}s^{1-N}U(s)^{-2}
\left(\int_0^s \tau^{N-1}U(\tau)^2d\tau\right)ds. 
\label{eq:3.5}
\end{eqnarray}
By \eqref{eq:2.10} we have 
\begin{equation}
\label{eq:3.6}
F(r)\preceq U(r)\int_0^r s^{1-N}(1+s)^{2A}
\left(\int_0^s \tau^{N-1}(1+\tau)^{-2A}d\tau\right)ds
\preceq r^2U(r)
\end{equation}
for all $r\ge 0$. 
Furthermore, since $F=F(|x|)$ satisfies  
$$
\Delta F-V(|x|)F=U(|x|)\quad\mbox{in}\quad{\bf R}^N,
$$
by \eqref{eq:3.4} we have 
\begin{eqnarray}
 & & \partial_t v-\Delta v+V(|x|)v
=\frac{A}{2}(T_2+t)^{\frac{A}{2}-1}\left[U(|x|)+\frac{A}{2}(T_2+t)^{-1}F(|x|)\right]\nonumber\\
 & & \qquad\qquad\qquad\qquad\qquad\quad
 -\frac{A}{2}(T_2+t)^{\frac{A}{2}-2}F(|x|)
 -\frac{A}{2}(T_2+t)^{\frac{A}{2}-1}U(|x|)\nonumber\\
 & & \qquad\qquad\qquad\qquad\quad
 =\left(\frac{A^2}{4}-\frac{A}{2}\right)(T_2+t)^{\frac{A}{2}-2}F(|x|)\ge 0
\label{eq:3.7}
\end{eqnarray}
for all $(x,t)\in{\bf R}^N\times(0,\infty)$. 
On the other hand, since $A<0$, 
by~\eqref{eq:1.3} we have $\omega>0$, 
and by condition~$(V)$~(ii) we can find a constant $L$ such that 
\begin{equation}
\label{eq:3.8}
V(r)>0
\end{equation}
for all $r\ge L$. 
Furthermore, taking a sufficiently large $L$ if necessary, 
by \eqref{eq:2.9} we have 
\begin{equation}
\label{eq:3.9}
0<\frac{|A|}{2}r^{-A-1}\le U'(r)\le 2|A|r^{-A-1}
\end{equation}
for all $r\ge L$. 
Then, similarly to \eqref{eq:3.6}, by~\eqref{eq:2.10} and \eqref{eq:3.9} 
we have 
\begin{eqnarray}
F'(r) \!\!\! & = &\!\!\!U'(r)\int_0^r s^{1-N}U(s)^{-2}
\left(\int_0^s \tau^{N-1}U(\tau)^2d\tau\right)ds\nonumber\\
 & & \qquad\qquad\qquad
+U(r)r^{1-N}U(r)^{-2}\int_0^r \tau^{N-1}U(\tau)^2d\tau\nonumber\\
\!\!\! & \preceq &\!\!\! 
r^2U'(r)+rU(r)\preceq r^{-A+1}\label{eq:3.10}
\end{eqnarray}
for all $r\ge L$. 

Let $\eta$ be a sufficiently small positive constant. 
Then, by \eqref{eq:3.6} we have 
\begin{equation}
\label{eq:3.11}
\frac{|A|}{2}(T_2+t)^{-1}F(|x|)
\le C_1(T_2+t)^{-1}|x|^2 U(|x|)
\le 4C_1\eta^2 U(|x|)\le\frac{1}{2}U(|x|)
\end{equation}
for all $(x,t)\in{\bf R}^N\times(0,\infty)$ with $|x|\le 2\eta(T_2+t)^{1/2}$, 
where $C_1$ is a constant. 
This together with \eqref{eq:2.10} and \eqref{eq:3.4} implies that 
\begin{eqnarray}
\frac{1}{2}(T_2+t)^{\frac{A}{2}}U(|x|) \!\!\! & \le &\!\!\! (1-4C_1\eta^2)(T_2+t)^{\frac{A}{2}}U(|x|)\notag\\
\!\!\! & \le &\!\!\! v(|x|,t)\notag\\
\!\!\! & \le &\!\!\! (1+4C_1\eta^2)(T_2+t)^{\frac{A}{2}}U(|x|)
\le\frac{3}{2}(T_2+t)^{\frac{A}{2}}U(|x|)\preceq 1
\label{eq:3.12}
\end{eqnarray}
for all $(x,t)\in{\bf R}^N\times(0,\infty)$ with $|x|\le 2\eta(T_2+t)^{1/2}$. 
In particular, taking sufficiently large $L$ and $T_2$ and a sufficiently small $\eta$ if necessary, 
by \eqref{eq:2.10} and \eqref{eq:3.12} we have 
\begin{eqnarray}
 & & v(|x|,t)=(1+O(\eta^2))(T_2+t)^{\frac{A}{2}}U(|x|)\nonumber\\
 & & \qquad\quad\,\,\,
 =(1+O(\eta^2))(T_2+t)^{\frac{A}{2}}|x|^{-A}(1+o(1))\nonumber\\
\label{eq:3.13}
 & & \qquad\quad\,\,\,
=(1+O(\eta^2))\eta^{-A}(1+o(1))
\le (5\eta/4)^{-A}
\end{eqnarray}
for all $(x,t)\in{\bf R}^N\times(0,\infty)$ with $|x|=\eta(T_2+t)^{1/2}\ge L$. 
Similarly, we have 
\begin{equation}
\label{eq:3.14}
v(|x|,t)=(1+O(\eta^2))(2\eta)^{-A}(1+o(1))
\ge (7\eta/4)^{-A}
\end{equation}
for all $(x,t)\in{\bf R}^N\times(0,\infty)$ with $|x|=2\eta(T_2+t)^{1/2}$.
Furthermore, 
taking a sufficiently small $\eta$ if necessary, 
by \eqref{eq:3.9} and \eqref{eq:3.10} we have 
$$
\frac{|A|}{2}(T_2+t)^{-1}F'(|x|)\le C_2(T_2+t)^{-1}|x|^{-A+1}
\le 4C_2\eta^2 |x|^{-A-1}<U'(|x|)
$$
for all $(x,t)\in{\bf R}^N\times(0,\infty)$ with $\eta(T_2+t)^{1/2}\le |x|\le 2\eta(T_2+t)^{1/2}$, 
where $C_2$ is a constant. 
This implies that 
\begin{equation}
\label{eq:3.15}
\frac{\partial}{\partial r}v(|x|,t)
=(T_2+t)^{\frac{A}{2}}
\left[U'(|x|)+\frac{A}{2}(T+t)^{-1}F'(|x|)\right]>0
\end{equation}
for all $(x,t)\in{\bf R}^N\times(0,\infty)$ with $\eta(T_2+t)^{1/2}\le |x|\le 2\eta(T_2+t)^{1/2}$. 
By \eqref{eq:3.13}--\eqref{eq:3.15}
we can take a smooth function $\rho=\rho(t)$ on $[0,\infty)$ such that 
$$
\eta(T_2+t)^{1/2}\le\rho(t)\le 2\eta(T_2+t)^{1/2}
\qquad\mbox{and}\qquad
v(\rho(t),t)=(3\eta/2)^{-A}
$$
for all $t>0$. 

For any sufficiently large $\kappa>0$, we define 
\begin{equation}
\label{eq:3.16}
W_2(x,t):=
\left\{
\begin{array}{ll}
\kappa v(x,t) & \mbox{if}\quad |x|\le\rho(t),\vspace{5pt}\\
\kappa (3\eta/2)^{-A} & \mbox{if}\quad |x|>\rho(t),
\end{array}
\right.
\end{equation}
for $(x,t)\in{\bf R}^N\times(0,\infty)$. 
Then, by \eqref{eq:2.10}, \eqref{eq:3.12}, and \eqref{eq:3.16} we have 
\begin{equation}
\label{eq:3.17}
W_2(x,0)\ge 1\quad\mbox{in}\quad{\bf R}^N. 
\end{equation}
Furthermore, 
since $\rho(t)\ge L$, by \eqref{eq:3.8} and \eqref{eq:3.16} we have 
\begin{equation}
\label{eq:3.18}
\partial_t W_2-\Delta W_2+V(|x|)W_2
=V(|x|)\kappa(3\eta/2)^{-A}>0
\end{equation}
for all $(x,t)\in{\bf R}^N\times(0,\infty)$ with $|x|>\rho(t)$. 
Then, by \eqref{eq:3.7}, \eqref{eq:3.15}, and \eqref{eq:3.18}
we see that the function $W_2$ is a supersolution of \eqref{eq:1.6}, 
and by \eqref{eq:3.17} we obtain 
\begin{equation}
\label{eq:3.19}
e^{-tH}1\le W_2(x,t),\qquad(x,t)\in{\bf R}^N\times(0,\infty). 
\end{equation}
In addition, by \eqref{eq:3.12} and \eqref{eq:3.16}, 
for any $\delta>0$, 
we have 
\begin{equation}
\label{eq:3.20}
\chi_\delta(x,t)W_2(x,t)\preceq 1
\end{equation}
for all $(x,t)\in{\bf R}^N\times(0,\infty)$. 
Therefore, by \eqref{eq:3.19} and \eqref{eq:3.20} 
we apply the same argument as in the case $A>0$, and obtain \eqref{eq:3.1} in the case $A<0$. 
Thus Lemma~\ref{Lemma:3.1} follows. 
$\Box$\vspace{5pt}%

Next we give the decay estimates of $[e^{-tH}\phi](x)$ outside parabolic cones 
by using Lemma~\ref{Lemma:3.1} and the $L^\infty_{loc}$ estimates for 
parabolic equations. 
\begin{lemma}
\label{Lemma:3.2}
Assume the same conditions as in Theorem~{\rm\ref{Theorem:1.1}}. 
Let $(p,p,\sigma,\sigma)\in\Lambda$ with $2<p\le \infty$. 
Then, for any sufficiently small $\delta>0$, 
there exists a constant $C$ such that 
\begin{equation}
\label{eq:3.21}
\left|[e^{-tH}\phi](x)\right|\le Ct^{-\frac{N}{2p}}\|\phi\|_{L^{p,\sigma}},
\qquad
\phi\in L^{p,\sigma}({\bf R}^N),
\end{equation}
for all $(x,t)\in{\bf R}^N\times[2,\infty)$ with $|x|\ge\delta(1+t)^{1/2}$. 
\end{lemma}
{\bf Proof.} 
Let $\delta$ be a sufficiently small positive constant. 
Define 
$$
u(x,t):=[e^{-tH}\phi](x),\qquad
h(t):=2\delta(1+t)^{1/2}.
$$
Then we may
assume, without loss of generality, that 
\begin{equation}
\label{eq:3.22}
h(t)^2<\frac{t}{2}\quad\mbox{for}\quad t\ge 2.
\end{equation}
Let $(x_0,t_0)\in{\bf R}^N\times[2,\infty)$ with $|x_0|\ge h(t_0)$.
Then the function 
\begin{equation}
\label{eq:3.23}
\tilde{u}(z,\tau):=u(\eta z+x_0,\eta^2\tau+t_0)
\quad\mbox{with}\quad\eta=\frac{h(t_0)}{2}
\end{equation}
is defined in $B(0,1)\times(-1,0)$ and it satisfies 
\begin{equation}
\label{eq:3.24}
\partial_{\tau} \tilde{u}=\Delta\tilde{u}-\eta^2V(|\eta z+x_0|)\tilde{u}
 \quad\mbox{in}\quad B(0,1)\times(-1,0).
\end{equation}
Since 
\begin{equation}
\label{eq:3.25}
|\eta z+x_0|\ge |x_0|-\eta
\ge h(t_0)-\eta=\frac{1}{2}h(t_0)=\eta,
\qquad z\in B(0,1),
\end{equation}
by condition~$(V)$ we have 
\begin{equation}
\label{eq:3.26}
\left|\eta^2V(|\eta z+x|)\right|\le\frac{C_1\eta^2}{|\eta z+x|^2}
\le C_1,\qquad z\in B(0,1),
\end{equation}
for some constant $C_1$. 
Then, by \eqref{eq:3.24} and \eqref{eq:3.26}  
we apply the standard $L^\infty_{loc}$ estimates for 
parabolic equations 
to obtain 
$$
|\tilde{u}(0,0)|
\le C_2\sup_{-1<\tau<0}\|\tilde{u}(\tau)\|_{L^1(B(0,1))}
$$
for some constant $C_2$. 
This together with \eqref{eq:3.22}, \eqref{eq:3.23}, and the H\"older inequality in the Lorentz spaces implies 
\begin{eqnarray}
|u(x_0,t_0)| \!\!\! & = &\!\!\! |\tilde{u}(0,0)|
\le C_2\eta^{-N}\sup_{t_0-\eta^2<s<t_0}\|u(s)\|_{L^1(B(x_0,\eta))}\nonumber\\
\!\!\! & \le &\!\!\! C_2\eta^{-N}\|\chi_{B(x_0,\eta)}\|_{L^{p',\sigma'}}
\sup_{t_0/2<s<t_0}\|\chi_{B(x_0,\eta)}u(s)\|_{L^{p,\sigma}}.
\label{eq:3.27}
\end{eqnarray}
On the other hand, by \eqref{eq:2.1} we have 
\begin{equation}
\label{eq:3.28}
\|\chi_{B(x_0,\eta)}\|_{L^{p',\sigma'}}
=\left(\int_{B(0,\eta)}|x|^{\frac{N\sigma'}{p'}-N}dx\right)^{1/\sigma'}
\le C_3\eta^{\frac{N}{p'}}=C_3\eta^{N-\frac{N}{p}}
\end{equation}
for some constant $C_3$. 
Furthermore, by \eqref{eq:3.25} we have 
\begin{equation}
\label{eq:3.29}
B(x_0,\eta)\subset\left\{x\in{\bf R}^N\,:\,|x|\ge\frac{1}{2}h(t_0)\right\}
\subset
\left\{x\in{\bf R}^N\,:\,|x|\ge\delta(1+s)^{1/2}\right\}
\end{equation}
for all $t_0/2<s<t_0$. 
Therefore, by Lemma~\ref{Lemma:3.1} and \eqref{eq:3.27}--\eqref{eq:3.29} 
we have 
\begin{equation}
\label{eq:3.30}
|u(x_0,t_0)|
\preceq \eta^{-\frac{N}{p}}\sup_{t_0/2<s<t_0}\|\chi_\delta(s)u(s)\|_{L^{p,\sigma}}
\preceq t_0^{-\frac{N}{2p}}\|\phi\|_{L^{p,\sigma}}
\end{equation}
for all $(x_0,t_0)\in{\bf R}^N\times[2,\infty)$ 
with $|x_0|\ge h(t_0)\ge\delta(1+t_0)^{1/2}$. 
Thus we have \eqref{eq:3.21}, and 
the proof
 is complete. 
$\Box$
\vspace{3pt}

Now we are ready to complete the proof of Proposition~\ref{Proposition:3.1}. 
\vspace{5pt}
\newline
{\bf Proof of Proposition~\ref{Proposition:3.1}.} 
Let $(p,p,\sigma,\sigma)\in\Lambda$ and $2<p\le\infty$. 
Let $\phi\in L^{p,\sigma}({\bf R}^N)$ with $\|\phi\|_{L^{p,\sigma}}=1$, 
and set 
$u(x,t):=[e^{-tH}\phi](x)$. 
For any sufficiently small $\delta>0$, 
by Lemma~\ref{Lemma:3.1} we have
\begin{equation}
\label{eq:3.31}
\|u(t)\|_{L^{p,\sigma}}\preceq\|(1-\chi_\delta(t))u(t)\|_{L^{p,\sigma}}+\|\chi_\delta(t)u(t)\|_{L^{p,\sigma}}
\preceq \|(1-\chi_\delta(t))u(t)\|_{L^{p,\sigma}}+1
\end{equation}
for all $t\ge 2$. 
On the other hand, 
by \eqref{eq:2.5} we have 
\begin{equation}
\label{eq:3.32}
\|u(2)\|_{L^\infty({\bf R}^N)}\le C_1
\end{equation}
for some constant $C_1$. 
By Lemma~\ref{Lemma:3.2} and \eqref{eq:3.32} 
we have 
$$
|u(x,t)|\le C_2(1+t)^{-\frac{N}{2p}},
\qquad (x,t)\in \Gamma_\delta(2),
$$
for some constant $C_2$. 
Let $w$ be the supersolution given in Proposition~\ref{Proposition:2.4} 
with $T=2$, $\epsilon=\delta$, $\gamma_1=-N/2p$, and $\gamma_2=0$. 
Then, applying the comparison principle, we obtain  
\begin{equation}
\label{eq:3.33}
|u(x,t)|\le C_3w(x,t)\preceq t^{-\frac{N}{2p}+\frac{A}{2}}U(|x|),
\qquad
(x,t)\in D_\delta(2),
\end{equation}
for some constant $C_3$. 
Therefore, if $A>0$  and $1\le \sigma <\infty$, then, 
by \eqref{eq:2.1}, \eqref{eq:2.10}, and \eqref{eq:3.33} we have 
\begin{eqnarray}
 & & \|(1-\chi_\delta(t))u(t)\|_{L^{p,\sigma}}
\preceq t^{-\frac{N}{2p}+\frac{A}{2}}\|(1-\chi_\delta(t))U\|_{L^{p,\sigma}}
\nonumber\\
 & & \qquad\quad
 \asymp 
 t^{-\frac{N}{2p}+\frac{A}{2}}
\left(\int_{|x|<\delta(1+t)^{1/2}}\frac{|x|^{\frac{\sigma N}{p}-N}}{(1+|x|)^{A\sigma}}dx\right)^{1/\sigma}
\nonumber\\
 & & \qquad\quad
 \preceq
 \left\{
 \begin{array}{ll}
 1 & \mbox{if}\quad 2<p<\beta,
 \vspace{3pt}\\
 (\log t)^{\frac{1}{\sigma}} & \mbox{if}\quad p=\beta,
 \vspace{5pt}\\
 t^{-\frac{N}{2p}+\frac{A}{2}} & \mbox{if}\quad \beta<p< \infty,
 \label{eq:3.34}
 \end{array}
 \right.
\end{eqnarray}
for all $t\ge 2$. 
Similarly, if $A>0$ and $\sigma=\infty$, then we have 
\begin{eqnarray}
 & & \|(1-\chi_\delta(t))u(t)\|_{L^{p,\infty}}
\preceq t^{-\frac{N}{2p}+\frac{A}{2}}\|(1-\chi_\delta(t))U\|_{L^{p,\infty}}
\nonumber\\
 & & \qquad\quad
 \asymp 
 t^{-\frac{N}{2p}+\frac{A}{2}}
\sup_{|x|<\delta(1+t)^{1/2}}\frac{|x|^{\frac{N}{p}}}{(1+|x|)^{A}}
\nonumber\\
 & & \qquad\quad
 \preceq
 \left\{
 \begin{array}{ll}
 1 & \mbox{if}\quad 2<p\le \beta,
 \vspace{3pt}\\
 t^{-\frac{N}{2p}+\frac{A}{2}} & \mbox{if}\quad \beta<p\le \infty,
 \label{eq:3.35}
 \end{array}
 \right.
\end{eqnarray}
for all $t\ge 2$. 
Therefore, by \eqref{eq:3.31}, \eqref{eq:3.34}, and \eqref{eq:3.35} we have assertion (II) for $2<p\le \infty$.

On the other hand, if $A\le 0$, then, by \eqref{eq:2.10} and \eqref{eq:3.33} we have 
$$
|u(x,t)|\preceq t^{-\frac{N}{2p}+\frac{A}{2}}(1+|x|)^{-A}\preceq t^{-\frac{N}{2p}},
\qquad
(x,t)\in D_\delta(2).
$$
This together with \eqref{eq:2.1} yields 
\begin{equation}
\label{eq:3.36}
\|(1-\chi_\delta(t))u(t)\|_{L^{p,\sigma}}\preceq
\left\{
\begin{array}{l}
 \displaystyle{t^{-\frac{N}{2p}}
\left(\int_{|x|<\delta(1+t)^{1/2}}|x|^{\frac{\sigma N}{p}-N}dx\right)^{1/\sigma}}\preceq 1, \quad 1\le \sigma <\infty,\vspace{5pt}\\
  \displaystyle{t^{-\frac{N}{2p}}
\sup_{|x|<\delta(1+t)^{1/2}}|x|^{\frac{N}{p}}}\preceq 1, \hspace{80pt}  \sigma=\infty,
\end{array}
\right.
\end{equation}
for all $t\ge 2$. 
Therefore, by \eqref{eq:3.31} and \eqref{eq:3.36}
we have assertions~(I) for $2<p\le \infty$. 
Finally, by Proposition \ref{Proposition:2.1}~(i) and Proposition \ref{Proposition:2.2} 
we have assertions~(I) and (II) for $1\le p\le\infty$, and 
the proof of Proposition~\ref{Proposition:3.1} is complete. 
$\Box$
\section{Decay estimates of $\|e^{-tH}\|_{(L^{p,\sigma}\to L^{q,\theta})}$} 
In this section we prove the following proposition on the decay rates of $\|e^{-tH}\|_{(L^{p,\sigma}\to L^{q,\theta})}$, 
which is a generalization of \cite[Proposition~4.1]{IIY}.
\begin{proposition}
\label{Proposition:4.1}
Assume the same conditions as in Theorem~{\rm\ref{Theorem:1.1}}. 
Let 
$$
{\rm(i)}\quad 1\le p\le r< q\le\infty,\,\,\,\,\sigma,\theta\in[1,\infty]
\,\,\quad\mbox{or}\,\,\quad
{\rm (ii)}\quad1\le p\le r=q\le \infty,\,\,\,\, 1\le \sigma\le \theta \le \infty.
$$ 
Assume that there exist constants $d_1$ and $d_2$ such that 
\begin{equation}
\label{eq:4.1}
\|e^{-tH}\|_{(L^{p,\sigma}\to L^{r,\sigma})} \preceq t^{d_1}(\log t)^{d_2}
\end{equation}
for all $t\ge 2$. 
Then 
\begin{equation}
\label{eq:4.2}
\|e^{-tH}\|_{(L^{p,\sigma}\to L^{q,\theta})}\preceq
t^{d_1}(\log t)^{d_2}\times 
\left\{
\begin{array}{ll}
t^{-\frac{N}{2}(\frac{1}{r}-\frac{1}{q})} & \mbox{if}\quad\mbox{$q<\infty$ and $Aq<N$},\vspace{3pt}\\
t^{-\frac{N}{2r}} & \mbox{if}\quad\mbox{$q=\infty$ and $A\le 0$},\vspace{3pt}\\
t^{-\frac{N}{2r}+\frac{A}{2}}(\log t)^{\frac{1}{\theta}}& \mbox{if}\quad \mbox{$q<\infty$ and $Aq=N$},\vspace{3pt}\\
t^{-\frac{N}{2r}+\frac{A}{2}} & \mbox{otherwise},
\end{array}
\right.
\end{equation}
for all $t\ge 2$. 
\end{proposition}
In order to prove Proposition~\ref{Proposition:4.1}, 
we prepare the following lemma. 
\begin{lemma}
\label{Lemma:4.1}
Assume the same conditions as in Proposition~{\rm\ref{Proposition:4.1}}. 
Then, for any sufficiently small $\delta>0$, 
\begin{equation}
\label{eq:4.3}
\left|[e^{-tH}\phi](x)\right|\preceq t^{d_1-\frac{N}{2r}}(\log t)^{d_2}\|\phi\|_{L^{p,\sigma}},
\qquad
\phi\in L^{p,\sigma}({\bf R}^N),
\end{equation}
for all $(x,t)\in{\bf R}^N\times[2,\infty)$ with $|x|\ge\delta(1+t)^{1/2}$. 
Furthermore, 
\begin{equation}
\label{eq:4.4}
|[e^{-tH}\phi](x)|\preceq t^{d_1-\frac{N}{2r}+\frac{A}{2}}(\log t)^{d_2}\|\phi\|_{L^{p,\sigma}}U(|x|),
\qquad
\phi\in L^{p,\sigma}({\bf R}^N),
\end{equation}
for all $(x,t)\in D_\delta(2)$. 
\end{lemma}
{\bf Proof.}
This lemma is proved by a similar argument as in Section~3. 
Let $\delta$ be a sufficiently small positive constant. 
By the same argument as in \eqref{eq:3.30} with the aid of \eqref{eq:4.1}
we have
$$
 |u(x,t)|
 \preceq t^{-\frac{N}{2r}}\sup_{t/2<s<t}\|\chi_\delta(s)u(s)\|_{L^{r,\sigma}}
\preceq t^{d_1-\frac{N}{2r}}(\log t)^{d_2}\|\phi\|_{L^{p,\sigma}}
$$
for all $(x,t)\in{\bf R}^N\times[2,\infty)$ with $|x|\ge \delta(1+t)^{1/2}$, 
and obtain \eqref{eq:4.3}. 
Furthermore, similarly to \eqref{eq:3.33}, 
by Proposition~\ref{Proposition:2.4} with $\gamma_1=d_1-N/2r$ and $\gamma_2=d_2$ 
we apply the comparison principle to obtain  
$$
|u(x,t)|\preceq t^{d_1-\frac{N}{2r}+\frac{A}{2}}(\log t)^{d_2}\|\phi\|_{L^{p,\sigma}}U(|x|)
$$
for all $(x,t)\in D_\delta(2)$. This implies \eqref{eq:4.4}. 
Thus Lemma~\ref{Lemma:4.1} follows.  
$\Box$\vspace{5pt}
\newline
{\bf Proof of Proposition~\ref{Proposition:4.1}.} 
Let $(p,q,\sigma,\theta)\in\Lambda$ be such that $1\le p\le r\le q< \infty$ and $1\le\theta<\infty$. 
Let $\phi\in L^{p,\sigma}({\bf R}^N)$ with $\|\phi\|_{L^{p,\sigma}}=1$ and 
set 
$u(t):=e^{-tH}\phi$. 
Assume \eqref{eq:4.1}. 
Then, for any $\delta>0$, by \eqref{eq:2.1} we have
\begin{eqnarray}
 & & \|u(t)\|_{L^{q,\theta}}
 \le 
 \left(
  \int_{|x|<\delta(1+t)^{1/2}}\left(|x|^{N/q}u^{\sharp}\right)^{\theta}\frac{dx}{|x|^N}
 \right)^{1/\theta}\notag\\
 & & \qquad\qquad\qquad\qquad
 +\left(
  \int_{|x|\ge \delta(1+t)^{1/2}}\left(|x|^{N/q}u^{\sharp}\right)^{\theta}\frac{dx}{|x|^N}
 \right)^{1/\theta}.
 \label{eq:4.5}
\end{eqnarray}
We first consider the case $A>0$. 
Since 
$$
U(|x|)\asymp (1+|x|)^{-A},\qquad x\in{\bf R}^N,
$$
Lemma~\ref{Lemma:4.1} implies 
$$
 u(t)^{\sharp}(x)
 \preceq 
  \left\{
   \begin{array}{l}
    t^{d_1-\frac{N}{2r}}(\log r)^{d_2}\vspace{3pt}
    \\
    \qquad\qquad
    \mbox{for all $(x,t)\in {\bf R}^N\times [2,\infty)$ with $|x|\ge \delta(1+t)^{1/2}$},\vspace{5pt}
    \\
    t^{d_1-\frac{N}{2r}+\frac{A}{2}}(\log r)^{d_2}(1+|x|)^{-A}\quad 
    \mbox{for all $(x,t)\in D_{\delta}(2)$}.
   \end{array}
  \right.
$$
Then, taking a sufficiently small $\delta$ if necessary,  
we have 
\begin{eqnarray}
 & & 
 \left(
  \int_{|x|<\delta(1+t)^{1/2}}\left(|x|^{N/q}u^{\sharp}\right)^{\theta}\frac{dx}{|x|^N}
 \right)^{1/\theta} \nonumber\\
 & & \qquad
\preceq 
t^{d_1-\frac{N}{2r}+\frac{A}{2}}(\log t)^{d_2}
\left(\int_{|x|<\delta(1+t)^{1/2}}\frac{|x|^{\theta N/q-N}}{(1+|x|)^{A\theta}}dx\right)^{1/\theta}\nonumber\\
 & & \qquad
 \preceq t^{d_1-\frac{N}{2r}+\frac{A}{2}}(\log t)^{d_2}\times
 \left\{
 \begin{array}{ll} 
 t^{\frac{N}{2q}-\frac{A}{2}} & \mbox{if}\quad q<\infty\quad\mbox{and}\quad Aq<N,\\
 (\log t)^{\frac{1}{\theta}} & \mbox{if}\quad q<\infty\quad\mbox{and}\quad Aq=N,\vspace{3pt}\\
 1 & \mbox{otherwise}, 
 \end{array}
 \right.
\label{eq:4.6}
\end{eqnarray}
for all $t\ge 2$. 
Furthermore,
if $q>r$, then, by \eqref{eq:2.1} and \eqref{eq:4.1} 
we see that  
\begin{eqnarray}
 & & \notag
 \left(
  \int_{|x|\ge \delta(1+t)^{1/2}}\left(|x|^{N/q}u^{\sharp}\right)^{\theta}\frac{dx}{|x|^N}
 \right)^{1/\theta}
 \\
 & & \notag
\qquad
\preceq 
 \left(
  \int_{|x|\ge\delta(1+t)^{1/2}}\left(|x|^{N/q}|x|^{-N/r}\|u(t)\|_{L^{r,\sigma}}\right)^{\theta}\frac{dx}{|x|^N}
 \right)^{1/\theta}
 \\
 & & \notag
\qquad
\preceq 
t^{d_1}(\log t)^{d_2}
 \left(
  \int_{|x|\ge \delta(1+t)^{1/2}}\left(|x|^{N/q-N/r}\right)^{\theta}\frac{dx}{|x|^N}
 \right)^{1/\theta}\\
 & &  \label{eq:4.7}
\qquad
\asymp
t^{d_1-\frac{N}{2}(\frac{1}{r}-\frac{1}{q})}(\log t)^{d_2}
\end{eqnarray}
for all $t\ge 2$. Therefore, by \eqref{eq:4.5}, \eqref{eq:4.6}, and \eqref{eq:4.7} we obtain 
$$
\|u(t)\|_{L^{q,\theta}} \preceq t^{d_1}(\log t)^{d_2}\times
\left\{
\begin{array}{ll}
t^{-\frac{N}{2}(\frac{1}{r}-\frac{1}{q})} & \mbox{if}\quad q<\infty\quad\mbox{and}\quad Aq<N,\\
t^{-\frac{N}{2r}+\frac{A}{2}}(\log(1+t))^{\frac{1}{\theta}} & \mbox{if}\quad q<\infty\quad\mbox{and}\quad Aq=N,\\
t^{-\frac{N}{2r}+\frac{A}{2}} & \mbox{otherwise},
\end{array}
\right.
$$
for all $t\ge 2$. This implies \eqref{eq:4.2} in the case~(i) with $1\le\theta<\infty$. 
On the other hand,
if $1\le p\le r=q< \infty$ and $1\le\sigma\le\theta<\infty$,
then, by \eqref{eq:4.1} we have
$$
 \left(
  \int_{|x|\ge \delta(1+t)^{1/2}}\left(|x|^{N/q}u^{\sharp}\right)^{\theta}\frac{dx}{|x|^N}
 \right)^{1/\theta}
\le 
\|u\|_{L^{q,\theta}}
\preceq
\|u\|_{L^{q,\sigma}}
\preceq 
t^{d_1}(\log t)^{d_2}
$$
for all $t\ge 2$. 
This together with \eqref{eq:4.5} and \eqref{eq:4.6} implies \eqref{eq:4.2} in the case~(ii) with $1\le\theta<\infty$. 
Therefore Proposition~\ref{Proposition:4.1} follows in the case where $A>0$ and $1\le\theta<\infty$. 

Next we consider the case $A\le 0$. 
By \eqref{eq:2.10}, \eqref{eq:4.3}, and \eqref{eq:4.4} we have
$$
 u(t)^{\sharp}(x)
 \preceq 
    t^{d_1-\frac{N}{2r}}(\log r)^{d_2}\|\phi \|_{L^{p,\sigma}}\quad 
    \mbox{for all $(x,t)\in {\bf R}^N\times [2,\infty)$}.
$$
This yields  
\begin{eqnarray*}
 & & 
 \left(
  \int_{|x|<\delta(1+t)^{1/2}}\left(|x|^{N/q}u^{\sharp}\right)^{\theta}\frac{dx}{|x|^N}
 \right)^{1/\theta} \nonumber\\
 & & 
\preceq 
t^{d_1-\frac{N}{2r}}(\log t)^{d_2}
\left(\int_{|x|<\delta(1+t)^{1/2}}|x|^{\theta N/q-N}dx\right)^{1/\theta}
\preceq 
t^{d_1-\frac{N}{2}(\frac{1}{r}-\frac{1}{q})}(\log t)^{d_2}
\end{eqnarray*}
for all $t\ge 2$. Then, by the same argument as in the case $A>0$ 
we have \eqref{eq:4.2}, and see that  
Proposition~\ref{Proposition:4.1} holds for the case where $A\le 0$ and $1\le\theta<\infty$. 
Thus Proposition~\ref{Proposition:4.1} follows for the case $1\le\theta<\infty$. 
Similarly, we can prove Proposition~\ref{Proposition:4.1} for the case $\theta=\infty$, 
and Proposition~\ref{Proposition:4.1} follows. 
$\Box$ 
\section{Proof of Theorem~\ref{Theorem:1.1}}
In this section we complete the proof of Theorem~\ref{Theorem:1.1}. 
We prepare the following proposition, 
which is useful to obtain the lower decay estimates of $\|e^{-tH}\|_{(L^{p,\sigma}\to L^{p,\sigma})}$. 
\begin{proposition}
\label{Proposition:5.1}
Assume the same conditions as in Theorem~{\rm\ref{Theorem:1.1}}. 
Let $\phi$ be a radially symmetric function in ${\bf R}^N$ such that $\phi\in C_0({\bf R}^N)$ 
and $\phi\ge(\not\equiv)\,0$ in ${\bf R}^N$. 
Then 
\begin{equation}
\label{eq:5.1}
\|e^{-tH}\phi\|_{L^1}\preceq t^{\frac{A}{2}},\qquad t\ge 2.
\end{equation}
Furthermore, for any sufficiently small $\epsilon>0$, 
\begin{equation}
\label{eq:5.2}
(e^{-tH}\phi)(x)\succeq t^{-\frac{N}{2}+A}U(|x|)
\end{equation}
for all $(x,t)\in{\bf R}^N\times[2,\infty)$ with $|x|\le\epsilon(1+t)^{1/2}$. 
\end{proposition}
{\bf Proof.}
We prove Proposition~\ref{Proposition:5.1} by a similar argument as in \cite{IK05}. 
Define
$$
u(x,t):=[e^{-tH}\phi](x),\quad
v(y,s):=(1+t)^{\frac{N}{2}}u(x,t),\quad
y:=(1+t)^{-\frac{1}{2}}x,\quad
s:=\log(1+t).
$$
Since $\phi\ge(\not\equiv)\,0$ in ${\bf R}^N$, 
we have 
\begin{equation}
\label{eq:5.3}
u(x,t)>0\qquad\mbox{in}\quad{\bf R}^N\times(0,\infty). 
\end{equation}
On the other hand, by \eqref{eq:2.7} we have
$$
\frac{d}{dt}\int_{{\bf R}^N}u(x,t)U(|x|)dx=0,\qquad t>0, 
$$
and obtain 
\begin{equation}
\label{eq:5.4}
\int_{{\bf R}^N}u(x,t)U(|x|)dx=\int_{{\bf R}^N}\phi(x)U(|x|)dx>0,\qquad t>0. 
\end{equation}
This together with \eqref{eq:2.10} implies 
\begin{eqnarray}
 & & \int_{{\bf R}^N}\phi(x)U(|x|)dx\ge \int_{(1+t)^{1/2}\le |x|\le 2(1+t)^{1/2}}u(x,t)U(|x|)dx\notag\\
 & & \qquad\quad
 \succeq (1+t)^{-\frac{A}{2}}\int_{(1+t)^{1/2}\le |x|\le 2(1+t)^{1/2}}u(x,t)dx\notag\\
 & & \qquad\quad
 \succeq (1+t)^{\frac{N}{2}-\frac{A}{2}}\min_{(1+t)^{1/2}\le |x|\le 2(1+t)^{1/2}}u(x,t)
\label{eq:5.5}
\end{eqnarray}
for all $t\ge e-1$. 
On the other hand, 
$v$ satisfies 
\begin{equation}
\label{eq:5.6}
\left\{
\begin{array}{ll}
\partial_s v=\Delta v+\displaystyle{\frac{y}{2}}\cdot\nabla v+\displaystyle{\frac{N}{2}}v-\tilde{V}(y,s)v & \vspace{3pt}\\
 \qquad
 =\displaystyle{\frac{1}{\rho}}\mbox{div}\,(\rho\nabla v)+\displaystyle{\frac{N}{2}}v-\tilde{V}(y,s)v
 & \quad\mbox{in}\quad{\bf R}^N\times(0,\infty),\vspace{3pt}\\
v(y,0)=\phi(y) & \quad\mbox{in}\quad{\bf R}^N,
\end{array}
\right.
\end{equation}
where $\rho(y)=e^{|y|^2/4}$ and $\tilde{V}(y,s)=e^s V(e^{s/2}y)$. 
Here, by condition~$(V)$ we have 
\begin{equation}
\label{eq:5.7}
|\tilde{V}(y,s)|\preceq |y|^{-2},
\qquad
\left|\tilde{V}(y,s)-\omega|y|^{-2}\right|
\preceq\frac{e^s}{(e^{s/2}|y|)^{2+a}}
\preceq e^{-\frac{a}{2}s}|y|^{-(2+a)},
\end{equation}
for all $(y,s)\in{\bf R}^N\times(0,\infty)$. 
Let $\epsilon$ be a sufficiently small positive constant. 
Then, by \eqref{eq:5.3} we apply 
the parabolic Harnack inequality to the solution $v$ of \eqref{eq:5.6}, 
and see that, for any $R\in(\epsilon,\infty)$, 
the inequality 
$$
\max_{\epsilon\le|y|\le R}v(y,s)\le C_1\min_{\epsilon\le|y|\le R} v(y,s+1),\qquad s\ge 2,
$$
holds for some positive constant $C_1$. 
This implies 
\begin{equation}
\label{eq:5.8}
\max_{\epsilon(1+t)^{1/2}\le|x|\le R(1+t)^{1/2}}u(x,t)\le C_2\min_{\epsilon(1+t')^{1/2}\le|x|\le R(1+t')^{1/2}}u(x,t')
\end{equation}
for all $t\ge e^2-1$, where $t'=e(1+t)-1$ and $C_2$ is a constant. 
Then, by \eqref{eq:5.5} and \eqref{eq:5.8}
we can find a positive constant $T_1$ such that 
\begin{equation}
\label{eq:5.9}
0\le u(x,t)\preceq (1+t')^{-\frac{N}{2}+\frac{A}{2}}\preceq (1+t)^{-\frac{N}{2}+\frac{A}{2}}
\end{equation}
for all $(x,t)\in{\bf R}^N\times(T_1,\infty)$ with $\epsilon(1+t)^{1/2}\le |x|\le 2(1+t)^{1/2}$. 
On the other hand, 
by \eqref{eq:2.5} we have 
\begin{equation}
\label{eq:5.10}
\|u(T_1)\|_{L^\infty({\bf R}^N)}\le C_3
\end{equation}
for some constant $C_3$. 
Let $w$ be the supersolution given in Proposition~\ref{Proposition:2.4} with $T=T_1$, $\gamma_1=-N/2+A/2$, and $\gamma_2=0$. 
By \eqref{eq:5.9} and \eqref{eq:5.10} 
we apply the comparison principle to obtain 
\begin{equation}
\label{eq:5.11}
0\le u(x,t)\preceq w(x,t)\preceq (1+t)^{-\frac{N}{2}+A}U(|x|)
\end{equation}
for all $(x,t)\in{\bf R}^N\times(T_1,\infty)$ with $|x|\le\epsilon(1+t)^{1/2}$. 
Furthermore, by \eqref{eq:5.7}  and \eqref{eq:5.11} 
we apply the same argument as in the proof of \cite[Lemma~4]{IK05} to the solution $v$ of \eqref{eq:5.6}, 
and obtain 
\begin{equation}
\label{eq:5.12}
\|v(s)\|_{L^2({\bf R}^N,\rho dy)}\preceq e^{\frac{A}{2}s}
\end{equation}
for all sufficiently large $s$. 
This implies 
\begin{equation}
\label{eq:5.13}
\|v(s)\|_{L^1({\bf R}^N)}
=\int_{{\bf R}^N}|v(s)|\rho^{1/2}\cdot \rho^{-1/2}dy
\preceq\left(\int_{{\bf R}^N}|v(s)|^2\rho dy\right)^{1/2}
\preceq e^{\frac{A}{2}s}
\end{equation}
for all sufficiently large $s$. 
This means that 
$$
\|u(t)\|_{L^1({\bf R}^N)}\preceq (1+t)^{\frac{A}{2}},
\qquad t\ge 2,
$$
and \eqref{eq:5.1} holds.
Furthermore, by \eqref{eq:2.10}, \eqref{eq:5.4}, and \eqref{eq:5.12}, 
taking a sufficiently small $\epsilon$ if necessary and applying the same argument as in \cite[Lemma~5]{IK05}, 
we can find positive constants $C_4$ and $L$ such that 
$$
\int_{\epsilon(1+t)^{1/2}\le|x|\le L(1+t)^{1/2}}u(x,t)U(|x|)dx\ge C_4>0
$$
for all sufficiently large $t$. 
This together with \eqref{eq:2.10} implies that 
\begin{equation}
\label{eq:5.14}
\max_{\epsilon(1+t)^{1/2}\le|x|\le L(1+t)^{1/2}}u(x,t)\succeq (1+t)^{-\frac{N}{2}+\frac{A}{2}}
\end{equation}
for all sufficiently large $t$. 
Therefore, by \eqref{eq:5.8} and \eqref{eq:5.14} we have 
\begin{equation}
\label{eq:5.15}
\min_{\epsilon(1+t)^{1/2}\le|x|\le L(1+t)^{1/2}}u(x,t)\succeq (1+t)^{-\frac{N}{2}+\frac{A}{2}},
\qquad t\ge T_2,
\end{equation}
for some constant $T_2$.

On the other hand, since $A<N/2$, 
the function 
$$
H(x,t):=(1+t)^{-\frac{N}{2}+A}U(|x|)
$$
satisfies 
\begin{equation}
\label{eq:5.16}
\partial_t H-\Delta H+VH=\left(-\frac{N}{2}+A\right)(1+t)^{-\frac{N}{2}+A-1}U<0\qquad\mbox{in}\quad{\bf R}^N\times(0,\infty).
\end{equation}
Furthermore, 
by \eqref{eq:2.10} we have 
\begin{equation}
\label{eq:5.17}
H(x,t)\preceq (1+t)^{-\frac{N}{2}+\frac{A}{2}}\qquad\mbox{on}\quad\Gamma_\epsilon(T_2).
\end{equation}
Then, by \eqref{eq:5.15}, \eqref{eq:5.16}, and \eqref{eq:5.17} 
we apply the comparison principle to obtain
$$
u(x,t)\succeq H(x,t)=(1+t)^{-\frac{N}{2}+A}U(|x|),
\qquad
(x,t)\in D_\epsilon(T_2).
$$
This together with \eqref{eq:5.3} implies \eqref{eq:5.2}. 
Thus Proposition~\ref{Proposition:5.1} follows.
$\Box$\vspace{5pt}

Now we are ready to complete the proof of Theorem~\ref{Theorem:1.1}. 
\vspace{3pt}
\newline
{\bf Proof of Theorem~\ref{Theorem:1.1}.}
Let $(p,q,\sigma,\theta)\in\Lambda$ and let $\phi$ be the function given in Proposition~\ref{Proposition:5.1}. 
We prove assertion~(I). 
Assume $A\le 0$. 
By Proposition~\ref{Proposition:3.1} we have \eqref{eq:4.1} with $d_1=d_2=0$ and $r=p$. 
Then it follows from Proposition~\ref{Proposition:4.1} that 
\begin{equation}
\label{eq:5.18}
\|e^{-tH}\|_{(L^{p,\sigma}\to L^{q,\theta})}\preceq t^{-\frac{N}{2}(\frac{1}{p}-\frac{1}{q})}
\end{equation}
for all $t\ge 2$. 
On the other hand, 
taking a sufficiently small $\epsilon>0$, 
by \eqref{eq:2.10} and \eqref{eq:5.2}
we have
$$
e^{-tH}\phi(x)\succeq t^{-\frac{N}{2}+\frac{A}{2}}
$$
for all $(x,t)\in{\bf R}^N\times[2,\infty)$ with $x\in E(t)$, 
where 
$$
E(t):=\{x\in{\bf R}^N\,:\,\epsilon(1+t)^{1/2}\le|x|\le 2\epsilon(1+t)^{1/2}\}.
$$ 
This implies that 
\begin{equation}
\label{eq:5.19}
\|e^{-tH}\phi\|_{L^{q,\theta}}\succeq t^{-\frac{N}{2}+\frac{A}{2}}\|\chi_{E(t)}\|_{L^{q,\theta}}
\succeq t^{-\frac{N}{2}(1-\frac{1}{q})+\frac{A}{2}},\qquad t\ge 2.  
\end{equation}
Furthermore, by \eqref{eq:5.1} and \eqref{eq:5.18} we have 
\begin{equation}
\label{eq:5.20}
\|e^{-2tH}\phi\|_{L^{p,\sigma}}\le\|e^{-tH}\|_{(L^1\to L^{p,\sigma})}\|e^{-tH}\phi\|_1
\preceq t^{-\frac{N}{2}(1-\frac{1}{p})+\frac{A}{2}}
\end{equation}
for all $t\ge 2$. 
Then, by \eqref{eq:5.18}--\eqref{eq:5.20} we have 
$$
t^{-\frac{N}{2}+\frac{A}{2}}\succeq
\|e^{-tH}\|_{(L^{p,\sigma}\to L^{q,\theta})}\ge\frac{\|e^{-3tH}\phi\|_{L^{q,\theta}}}{\|e^{-2tH}\phi\|_{L^{p,\sigma}}}
\succeq t^{-\frac{N}{2}(\frac{1}{p}-\frac{1}{q})}
$$
for all $t\ge 2$. 
Thus assertion~(I) follows. 
\vspace{3pt}

We prove assertion~(II).
Assume $A>0$. 
We first prove assertion~(II)~(i). 
Let $1\le p<\alpha$. 
By Proposition~\ref{Proposition:3.1}~(II) and \eqref{eq:2.6} 
we have 
$$
\|e^{-tH}\|_{(L^{p,\sigma}\to L^{p,\sigma})}
\preceq t^{-\frac{N}{2p'}+\frac{A}{2}}=t^{-\frac{N}{2}(1-\frac{1}{p})+\frac{A}{2}}
$$
for all $t\ge 2$. 
Then we apply Proposition~\ref{Proposition:4.1} with $r=p$ to obtain 
\begin{equation}
\label{eq:5.21}
\|e^{-tH}\|_{(L^{p,\sigma}\to L^{q,\theta})}\preceq
\left\{
\begin{array}{ll}
t^{-\frac{N}{2}(1-\frac{1}{q})+\frac{A}{2}} & \mbox{if}\quad\mbox{$p\le q<\beta$},\\
t^{-\frac{N}{2}+A}(\log t)^{\frac{1}{\theta}} & \mbox{if}\quad\mbox{$q=\beta$},\\
t^{-\frac{N}{2}+A} & \mbox{if}\quad\mbox{$\beta<q\le\infty$},
\end{array}
\right.
\end{equation}
for all $t\ge 2$. 
On the other hand, 
if $p\le q <\infty$ and $\theta<\infty$, then, by \eqref{eq:2.1}, \eqref{eq:2.10}, and \eqref{eq:5.2} 
we can find positive constants $T$ and $\epsilon$ such that 
\begin{eqnarray}
\|e^{-tH}\phi\|_{L^{q,\theta}} \!\!\!& \succeq &\!\!\! 
t^{-\frac{N}{2}+A}\left(\int_{|x|<\epsilon(1+t)^{1/2}}
 \left(|x|^{N/q}U^{\sharp}(x)\right)^{\theta}\frac{dx}{|x|^N}\right)^{1/\theta}\notag\\
 \!\!\!& \succeq &\!\!\! t^{-\frac{N}{2}+A}\left(\int_{|x|<\epsilon(1+t)^{1/2}}\frac{|x|^{N\theta/q-N}}{(1+|x|)^{A\theta}}dx\right)^{1/\theta}\notag\\
 \!\!\!& \succeq &\!\!\!
 \left\{
\begin{array}{ll}
t^{-\frac{N}{2}(1-\frac{1}{q})+\frac{A}{2}} &\mbox{if}\quad q<\beta,\\
t^{-\frac{N}{2}+A}(\log t)^{\frac{1}{\theta}} & \mbox{if}\quad q=\beta,\\
t^{-\frac{N}{2}+A} & \mbox{if}\quad \beta<q<\infty,
\end{array}
\right.\label{eq:5.22}
\end{eqnarray}
for all $t\ge T$. 
Similarly, if $p\le q \le \infty$ and $\theta=\infty$, then, by \eqref{eq:2.1}, \eqref{eq:2.10}, and \eqref{eq:5.2} we have
\begin{eqnarray}
\|e^{-tH}\phi\|_{L^{q,\infty}} \!\!\!& \succeq &\!\!\! 
t^{-\frac{N}{2}+A}
\sup_{|x|<\epsilon(1+t)^{1/2}}
 \left(
  |x|^{N/q}U^{\sharp}(x)
 \right)\notag\\
 \!\!\!& \succeq &\!\!\! t^{-\frac{N}{2}+A}\sup_{|x|<\epsilon(1+t)^{1/2}}\frac{|x|^{N/q}}{(1+|x|)^{A}}\notag\\
 \!\!\!& \succeq &\!\!\!
 \left\{
\begin{array}{ll}
t^{-\frac{N}{2}(1-\frac{1}{q})+\frac{A}{2}} &\mbox{if}\quad q<\beta,\\
t^{-\frac{N}{2}+A} & \mbox{if}\quad q=\beta,\\
t^{-\frac{N}{2}+A} & \mbox{if}\quad \beta<q\le\infty.
\end{array}
\right.\label{eq:5.23}
\end{eqnarray}
Therefore we deduce from \eqref{eq:5.22} and \eqref{eq:5.23} that 
\begin{equation}
\label{eq:5.24}
\|e^{-tH}\|_{(L^{p,\sigma}\to L^{q,\theta})}\ge\frac{\|e^{-tH}\phi\|_{L^{q,\theta}}}{\|\phi\|_{L^{p,\sigma}}}\succeq
\left\{
\begin{array}{ll}
t^{-\frac{N}{2}(1-\frac{1}{q})+\frac{A}{2}} & \mbox{if}\quad\mbox{$p\le q<\beta$},\\
t^{-\frac{N}{2}+A}(\log t)^{\frac{1}{\theta}} & \mbox{if}\quad\mbox{$q=\beta$},\\
t^{-\frac{N}{2}+A} & \mbox{if}\quad\mbox{$\beta<q\le\infty$},
\end{array}
\right.
\end{equation}
for all $t\ge 2$. Then assertion~(II)~(i) follows from \eqref{eq:5.21} and \eqref{eq:5.24}. 
Furthermore, assertion~(II)~(v) follows from assertion~(II)~(i) and Proposition~\ref{Proposition:2.2}. 
\vspace{3pt}

Next we prove assertion~(II)~(iii). 
Let $\alpha<p<\beta$. 
Similarly to \eqref{eq:5.21}, by Propositions~\ref{Proposition:3.1} and \ref{Proposition:4.1} 
we can obtain assertion~(II)~(iii) with $\asymp$ replaced by $\preceq$. 
Furthermore, 
by assertion~(II)~(i), \eqref{eq:5.22}, and \eqref{eq:5.23} 
we have 
\begin{eqnarray*}
 & & \|e^{-tH}\|_{(L^{p,\sigma}\to L^{q,\theta})}\ge\frac{\|e^{-2tH}\phi\|_{L^{q,\theta}}}{\|e^{-tH}\phi\|_{L^{p,\sigma}}}
\ge\frac{\|e^{-2tH}\phi\|_{L^{q,\theta}}}{\|e^{-tH}\|_{(L^1 \to L^{p,\sigma})}\|\phi\|_1}\\
 & & \qquad\qquad\,\,\,
\succeq t^{\frac{N}{2}(1-\frac{1}{p})-\frac{A}{2}}\frac{\|e^{-2tH}\phi\|_{L^{q,\theta}}}{\|\phi\|_1}
\succeq\left\{
\begin{array}{ll}
t^{-\frac{N}{2}(\frac{1}{p}-\frac{1}{q})} & \mbox{if}\quad p\le q<\beta,\\
t^{-\frac{N}{2p}+\frac{A}{2}}(\log t)^{\frac{1}{\theta}} &\mbox{if}\quad q=\beta,\\
t^{-\frac{N}{2p}+\frac{A}{2}} & \mbox{if}\quad \beta<q\le\infty,
\end{array}
\right.
\end{eqnarray*}
for all $t\ge 2$. Then assertion~(II)~(iii) follows. 
\vspace{3pt}

Next we prove assertions~(II)~(ii) and (iv). 
Let $p\in\{\alpha,\beta\}$. 
Due to Proposition~\ref{Proposition:2.2} and assertions~(II)~(i) and (iii), 
it suffices to consider the following three cases, 
$(p,q)=(\alpha,\alpha)$, $(\alpha,\beta)$, and $(\beta,\beta)$. 
In these three cases, by Propositions~\ref{Proposition:3.1} and \ref{Proposition:4.1}
we obtain the desired upper estimates of $\|e^{-tH}\|_{(L^{p,\sigma}\to L^{q,\theta})}$. 
It remains to prove 
\begin{eqnarray}
\label{eq:5.25}
 & & \|e^{-tH}\|_{(L^{\beta,\sigma}\to L^{\beta,\theta})}\asymp\|e^{-tH}\|_{(L^{\alpha,\theta'}\to L^{\alpha,\sigma'})}
 \succeq(\log t)^{\frac{1}{\theta}},\vspace{3pt}\\
 \label{eq:5.26}
 & & \|e^{-tH}\|_{(L^{\alpha,\sigma}\to L^{\beta,\theta})}\succeq t^{-\frac{N}{2}+A}(\log t)^{\frac{1}{\theta}+\frac{1}{\sigma'}},
\end{eqnarray}
for all $t\ge 2$. 

Let $\alpha<p<r<\beta$. 
It follows from Proposition~\ref{Proposition:2.1}~(ii)  that 
\begin{equation}
\label{eq:5.27}
\|e^{-tH}\|_{(L^{r,\sigma}\to L^{\beta,\theta})}\le\|e^{-tH}\|_{(L^{p,\sigma}\to L^{\beta,\theta})}^{1-\eta_1}
 \|e^{-tH}\|_{(L^{\beta,\sigma}\to L^{\beta,\theta})}^{\eta_1},
\end{equation}
where
$$
\frac{1}{r}=\frac{1-\eta_1}{p}+\frac{\eta_1}{\beta}. 
$$
On the other hand, by assertion~(II)~(iii) we have
\begin{equation}
\label{eq:5.28}
\|e^{-tH}\|_{(L^{r,\sigma}\to L^{\beta,\theta})}\asymp t^{-\frac{N}{2r}+\frac{A}{2}}(\log t)^{\frac{1}{\theta}},
\qquad
\|e^{-tH}\|_{(L^{p,\sigma}\to L^{\beta,\theta})}\asymp t^{-\frac{N}{2p}+\frac{A}{2}}(\log t)^{\frac{1}{\theta}}.
\end{equation}
Then, by \eqref{eq:5.27} and \eqref{eq:5.28} we obtain
$$
\|e^{-tH}\|_{(L^{\beta,\sigma}\to L^{\beta,\theta})}^{\eta_1}
\succeq t^{-\frac{N}{2r}+(1-\eta_1)\frac{N}{2p}+\frac{A}{2}\eta_1}(\log t)^{\frac{1}{\theta}\eta_1}
=t^{-\frac{N}{2\beta}\eta_1+\frac{A}{2}\eta_1}(\log t)^{\frac{\eta_1}{\theta}}=(\log t)^{\frac{\eta_1}{\theta}}
$$
for all $t\ge 2$. This together with Proposition~\ref{Proposition:2.2} yields \eqref{eq:5.25}.
\vspace{3pt}

We prove \eqref{eq:5.26}. 
The proof is divided into the following four cases: 
$$
(1)\ 1\le \sigma=\theta'\le \infty;
\ \ 
(2)\ 1\le\sigma<\theta'<\infty;
\ \ 
(3)\ 1<\theta'<\sigma\le\infty;
\ \ 
(4)\ 1\le \sigma \le \infty,\  \theta'\in\{1,\infty\}.
$$
We first consider the case~(1). 
By assertion~(II)~(ii) with $(p,q)=(\alpha,2)$ 
we have 
\begin{equation}
\label{eq:5.29}
\begin{split}
\|e^{-tH}\|_{(L^{\alpha,\theta'}\to L^{\beta,\theta})}
 &
 = 
\sup_{\|\phi\|_{L^{\alpha,\theta'}}=1}\|e^{-tH}\phi\|_{L^{\beta,\theta}}
\\
&
=\sup_{\|\phi\|_{L^{\alpha,\theta'}}=1}\sup_{\|\psi\|_{L^{\alpha,\theta'}}=1}\left|\int_{{\bf R}^N}[e^{-tH}\phi](x)\psi(x)dx\right|
\\
& 
\ge 
\sup_{\|\phi\|_{L^{\alpha,\theta'}}=1}\left|\int_{{\bf R}^N}[e^{-tH/2}\phi](x)[e^{-tH/2}\phi](x)dx\right|
\\
&
=\sup_{\|\phi\|_{L^{\alpha,\theta'}}=1}\|e^{-tH/2}\phi\|_2^2
\\
& 
= \|e^{-tH/2}\|_{(L^{\alpha,\theta'}\to L^2)}^2\asymp t^{-\frac{N}{2}+A}(\log t)^{\frac{2}{\theta}}
\end{split}
\end{equation}
for all $t\ge 2$. This implies \eqref{eq:5.26} in the case~(1). 

Next we consider the case~(2). 
It follows from Proposition~\ref{Proposition:2.1}~(iii) that 
\begin{equation}
\label{eq:5.30}
 \|e^{-tH}\|_{(L^{\alpha,\theta'}\to L^{\beta,\theta})}
  \le
    \|e^{-tH}\|^{1-\eta_2}_{(L^{\alpha,\sigma}\to L^{\beta,\theta})}
    \|e^{-tH}\|^{\eta_2}_{(L^{\alpha,\infty}\to L^{\beta,\theta})},
    \qquad t>0,
\end{equation}
where 
$$
\frac{1}{\theta'}=\frac{1-\eta_2}{\sigma}+\frac{\eta_2}{\infty}=\frac{1-\eta_2}{\sigma}. 
$$
By Propositions~\ref{Proposition:3.1} and \ref{Proposition:4.1} 
we have 
\begin{equation}
\label{eq:5.31}
\|e^{-tH}\|_{(L^{\alpha,\infty}\to L^{\beta,\theta})}
\preceq t^{-\frac{N}{2}+{A}}(\log t)^{\frac{1}{\theta}+1}
\end{equation}
for all $t\ge 2$.
Then, by \eqref{eq:5.29}, \eqref{eq:5.30}, and \eqref{eq:5.31}
we obtain 
$$
\|e^{-tH}\|_{(L^{\beta,\sigma}\to L^{\beta,\theta})}^{1-\eta_2}
\succeq 
t^{\left(-\frac{N}{2}+{A}\right)(1-\eta_2)}(\log t)^{\frac{2}{\theta}-\eta_2\left(\frac{1}{\theta}+1\right)}
=
\left(t^{-\frac{N}{2}+{A}}(\log t)^{\frac{1}{\theta}+\frac{1}{\sigma'}}\right)^{1-\eta_2}
$$
for all $t\ge 2$, which implies \eqref{eq:5.26} in the case~(2). 
Similarly, in the case~(3), 
we see that 
\begin{eqnarray*}
 & & \|e^{-tH}\|_{(L^{\alpha,\theta'}\to L^{\beta,\theta})}
 \le\|e^{-tH}\|^{1-\eta_3}_{(L^{\alpha,1}\to L^{\beta,\theta})}
 \|e^{-tH}\|^{\eta_3}_{(L^{\alpha,\sigma}\to L^{\beta,\theta})},\qquad t>0,\\
 & & \|e^{-tH}\|_{(L^{\alpha,1}\to L^{\beta,\theta})}
\preceq t^{-\frac{N}{2}+{A}}(\log t)^{\frac{1}{\theta}},
\qquad t\ge 2,
\end{eqnarray*}
where
$$
 \frac{1}{\theta'}=\frac{1-\eta_3}{1}+\frac{\eta_3}{\sigma}.
$$
These imply \eqref{eq:5.26} in the case~(3). 

Finally we consider the case~(4). 
If $1<\sigma<\infty$, then, by Proposition~\ref{Proposition:2.2} and \eqref{eq:5.26} in the cases~(2) and (3)
we have
\begin{eqnarray}
 & & \|e^{-tH}\|_{(L^{\alpha,\sigma}\to L^{\beta,\infty})}
\asymp 
\|e^{-tH}\|_{(L^{\alpha,1}\to L^{\beta,\sigma'})}
\succeq
t^{-\frac{N}{2}+A}(\log t)^{\frac{1}{\sigma'}},\label{eq:5.32}\\
 & & \|e^{-tH}\|_{(L^{\alpha,\sigma}\to L^{\beta,1})}
\asymp 
\|e^{-tH}\|_{(L^{\alpha,\infty}\to L^{\beta,\sigma'})}
\succeq 
t^{-\frac{N}{2}+A}(\log t)^{1+\frac{1}{\sigma'}},\label{eq:5.33}
\end{eqnarray}
for all $t\ge 2$. 
Furthermore, it follows from \eqref{eq:5.26} in the case~(1) that 
\begin{equation}
\label{eq:5.34}
\begin{split}
\|e^{-tH}\|_{(L^{\alpha,1}\to L^{\beta,\infty})}
&
\succeq t^{-\frac{N}{2}+A}, \\
\|e^{-tH}\|_{(L^{\alpha,\infty}\to L^{\beta,1})}
&
\succeq t^{-\frac{N}{2}+A}(\log t)^2,
\end{split}
\end{equation}
for all $t\ge 2$. 
It remains to prove the cases
$$
(\sigma,\theta')=(1,\infty)\qquad\mbox{and}\qquad
(\sigma,\theta')=(\infty,1).
$$
Let $1<\tilde{\sigma}<\infty$. 
Proposition~\ref{Proposition:2.1}~(iii) implies 
\begin{equation}
\label{eq:5.35}
\|e^{-tH}\|_{(L^{\alpha,\tilde{\sigma}}\to L^{\beta,\infty})} 
\le\|e^{-tH}\|^{1-\eta_4}_{(L^{\alpha,1}\to L^{\beta,\infty})}
\|e^{-tH}\|^{\eta_4}_{(L^{\alpha,\infty}\to L^{\beta,\infty})},
\end{equation}
where 
$$
\frac{1}{\tilde\sigma}=\frac{1-\eta_4}{1}+\frac{\eta_4}{\infty}=1-\eta_4.
$$
On the other hand, similarly to the case~(2), 
by Propositions~\ref{Proposition:3.1} and \ref{Proposition:4.1} 
we see that 
\begin{equation}
\label{eq:5.36}
\|e^{-tH}\|_{(L^{\alpha,1}\to L^{\beta,\infty})} \preceq t^{-\frac{N}{2}+A},\qquad t\ge 2.
\end{equation}
Furthermore, by \eqref{eq:5.32} we have
\begin{equation}
\label{eq:5.37}
\|e^{-tH}\|_{(L^{\alpha,\tilde{\sigma}}\to L^{\beta,\infty})} 
\succeq t^{-\frac{N}{2}+A}(\log t)^{\frac{1}{{\tilde{\sigma}}'}}=t^{-\frac{N}{2}+A}(\log t)^{\eta_4},
\qquad t\ge 2.
\end{equation}
Therefore we deduce from 
\eqref{eq:5.35}--\eqref{eq:5.37} and Proposition \ref{Proposition:2.2} that 
$$
\|e^{-tH}\|_{(L^{\alpha,1}\to L^{\beta,1})}
\asymp 
\|e^{-tH}\|_{(L^{\alpha,\infty}\to L^{\beta,\infty})}
\succeq t^{-\frac{N}{2}+A}\log t
$$
for all $t\ge 2$, 
which implies \eqref{eq:5.26} in the case~(4). 
Thus assertion~(II) follows, and the proof of Theorem~\ref{Theorem:1.1} is complete.
$\Box$

\end{document}